\documentclass[11pt,reqno]{amsart}
\usepackage[latin1]{inputenc} 
\usepackage{indentfirst}
\usepackage{epsfig} 
\usepackage{graphicx} 
\usepackage{multirow} 
\usepackage{verbatim} 
\usepackage{rotating} 
\usepackage{lettrine}
\usepackage{amscd}
\usepackage{pstricks}
\usepackage{pst-tree}
\usepackage{pst-node}
\usepackage{lscape}
\usepackage{ae}
\usepackage{float}
\usepackage[all]{xy}
\usepackage{icomma} 
\usepackage{varioref} 
\usepackage{layout}
\usepackage{enumerate}
\usepackage{amsfonts, amsmath, amssymb, stmaryrd, latexsym, dsfont, amsthm}
\usepackage{array}
\usepackage{longtable}
\usepackage{geometry}
\usepackage[frenchb,english]{babel}
\begin{document}
\renewcommand{\baselinestretch}{1.3}
\newcommand{\K}{\mathbb{Q}(\sqrt{d},i)}
\newcommand{\q}{\mathbb{Q}(\sqrt{d})/\mathbb{Q}}
\theoremstyle{plain}
\newtheorem*{mainthem}{Theorem}
\newtheorem{them}{Theorem}[]
\newtheorem{lem}{Lemma}[]
\newtheorem{propo}{Proposition}[]
\newtheorem{coro}{Corollary}[]
\newtheorem{proprs}{properties}[]
\theoremstyle{remark}
\newtheorem{rema}{Remark}[]
\newtheorem{remas}{Remarks}[]
\newtheorem{exam}{\textbf{Numerical Example}}[]
\newtheorem{exams}{\textbf{Numerical Examples}}[]
\newtheorem{df}{definition}[]
\newtheorem{dfs}{definitions}[]
\def\NN{\mathds{N}}
\def\RR{\mathbb{R}}
\def\HH{I\!\! H}
\def\QQ{\mathbb{Q}}
\def\CC{\mathbb{C}}
\def\ZZ{\mathbb{Z}}
\def\OO{\mathcal{O}}
\def\kk{\mathds{k}}
\def\KK{\mathbb{K}}
\def\ho{\mathcal{H}_0^{\frac{h(d)}{2}}}
\def\h{\overline{\mathcal{H}_0}}
\def\hh{\overline{\mathcal{H}_1}}
\def\hhh{\overline{\mathcal{H}_2}}
\def\jj{\overline{\mathcal{H}_3}}
\def\jjj{\overline{\mathcal{H}_4}}
\def\k{\mathds{k}^{(*)}}
\def\l{\mathds{L}}
\title[ On the strongly ambiguous...]{On the strongly ambiguous classes of $\kk/\QQ(i)$ where $\kk= \QQ(\sqrt{2p_1p_2},i)$}
\author[Abdelmalek AZIZI]{Abdelmalek Azizi}
\address{Abdelmalek Azizi et Abdelkader Zekhnini: Département de Mathématiques, Faculté des Sciences, Université Mohammed 1, Oujda, Morocco }
\author{Abdelkader Zekhnini}
\email{abdelmalekazizi@yahoo.fr}
\email{zekha1@yahoo.fr}
\author[Mohammed Taous]{Mohammed Taous}
\address{Mohammed Taous: Département de Mathématiques, Faculté des Sciences et Techniques, Université Moulay Ismail, Errachidia, Morocco}
\email{taousm@hotmail.com}

\maketitle
\selectlanguage{english}
\begin{abstract}
We construct an infinite family of imaginary bicyclic biquadratic number fields  $\kk$ with the 2-ranks of their 2-class groups are $\geq3$, whose  strongly ambiguous classes of $\kk/\QQ(i)$  capitulate in the absolute genus field $\k$, which is  strictly included in the relative genus field $(\kk/\QQ(i))^*$ and we study the capitulation  of the $2$-ideal classes of $\kk$ in its quadratic extensions included in $\k$.
\end{abstract}
 \section{\textbf{Introduction}.}
 Let $\mathds{k}$ be an algebraic number field and let $\mathbf{C}_{\mathds{k},2}$ denote its 2-class group, that is the 2-Sylow
subgroup of the ideal class group, $\mathbf{C}_\mathds{k}$,  of $\mathds{k}$.
We denote by  $\k$  the absolute genus field of $\mathds{k}$.
Suppose $\mathds{F}$ is a finite extension of $\mathds{k}$,
then we say that an ideal class of $\mathds{k}$ capitulates in $\mathds{F}$ if it is in the kernel of the
homomorphism
$$J_\mathds{F}: \mathbf{C}_\mathds{k} \longrightarrow \mathbf{C}_\mathds{F}$$ induced by extension of ideals from $\mathds{k}$
to $\mathds{F}$. An
important problem in Number Theory is to determine explicitly the kernel of  $J_\mathds{F}$, which is
usually called the capitulation kernel. The classical Principal Ideal Theorem  asserts that $kerJ_\mathds{F}$ is all $\mathbf{C}_\mathds{k}$ if $\mathds{F}$ is the Hilbert class field of $\kk$. If $\mathds{F}$ is
the  relative genus  field of a cyclic extension $\mathds{K}/\kk$, which we denote  by $(\mathds{K}/\kk)^*$ and that is  the maximal unramified extension  of  $\mathds{K}$ which is obtained by composing  $\mathds{K}$ and an abelian extension over $\kk$, F. Terada  states in  \cite{FT-71} that all the ambiguous ideal classes of $\mathds{K}/\kk$ capitulate in
 $(\mathds{K}/\kk)^*$; if $\mathds{F}$ is the  absolute genus  field of an abelian extension $\mathds{K}/\QQ$, then H. Furuya confirms in \cite{Fu-77} that every   strongly ambiguous  class, that is an ambiguous ideal class  represented by an ambiguous ideal, of $\mathds{K}/\QQ$ capitulate in  $\mathds{F}$. In this paper we construct a family of number field $\kk$ for which all the strongly ambiguous  classes of $\kk/\QQ(i)$ capitulate in  $\k\varsubsetneq (\kk/\QQ(i))^*$ and all classes that capitulate in an unramified quadratic extension $\mathds{K}$ of  $\kk$ that is abelian over $\QQ$ are strongly ambiguous  classes of $\kk/\QQ(i)$.\par
Let $p_1\equiv p_2\equiv1\pmod4$ be primes,  $\mathds{k}=\QQ(\sqrt{2p_1p_2},i)$ and
 $\mathds{K}$ be an unramified quadratic extension of $\kk$ that is abelian  over $\QQ$, let  $Am_s(\kk/\QQ(i))$
denote the group of the strongly ambiguous classes of  $\kk/\QQ(i)$. Our
 object  is to determine  the 2-classes of the field $\mathds{k}$ which capitulate in the extension
$\mathds{K}$ involving the fundamental units of the three quadratic  subfields of $\mathds{k}$, we prove that
$kerJ_\mathds{K}\subset Am_s(\kk/\QQ(i))$  and we infer that   $Am_s(\kk/\QQ(i))\subseteq kerJ_{\k}$. As
an application we will determine these 2-classes when  $\mathbf{C}_{\mathds{k},2}$
is of type $(2, 2, 2)$. This study is based on genus theory, class group theory and other theorems as the next result giving the number of classes which capitulate in a cyclic extension of prime degree: if $K/k$ is a cyclic extension of prime degree, then the number of classes which capitulate in $K/k$ is:
\begin{equation}\label{1} [K:k][E_k:N(E_K)],\end{equation}
 where $E_k$ and $E_K$ are the  unit groups of $k$ and $K$ respectively and $N$ is the norm of $K/k$.\par
   Let $m$ be a square-free integer and $K$ be a number field, during this paper, we adopt the following notations:
 \begin{itemize}
   \item   $p_1\equiv p_2\equiv1\pmod4$ are primes.
   \item $\mathds{k}$: denotes the field $\QQ(\sqrt{2p_1p_2},\sqrt{-1})$.
   \item $\mathcal{O}_K$: the ring of integers of $K$.
   \item $E_K$: the unit group of $\mathcal{O}_K$.
   \item  $W_K$: the group of roots of unity contained in $K$.
   \item $F.S.U$: the fundamental system of units.
   \item $K^+$: the maximal real subfield of $K$, if it is a CM-field.
   \item $Q_K=[E_K:W_KE_{K^+}]$ is Hasse's unit index, if $K$ is a CM-field.
   \item $q(K/\QQ)=[E_K:\prod_i^s E_{k_i}]$ is the unit index of $K$, if $K$ is multiquadratic, where $k_i$  are the quadratic subfields of $K$.
  \item $K^{(*)}$: the absolute genus field of $K$.
  \item $\mathbf{C}_{K,2}$: the 2-class group of $K$.
  \item $i=\sqrt{-1}$.
  \item $\varepsilon_m$: the fundamental unit of $\QQ(\sqrt m)$.
  \item $N(a)$: denotes the absolute norm of a number $a$.
 \end{itemize}
 Our main theorem is.
 \begin{mainthem}
Let $Am_s(\kk/\QQ(i))$
denote the group of the strongly ambiguous classes of  $\kk/\QQ(i)$.
If $\mathds{K}$ is an unramified quadratic extension of $\kk$ that is  abelian  over $\QQ$,
then
\begin{enumerate}[\upshape\indent(1)]
  \item $kerJ_\mathds{K}\subset Am_s(\kk/\QQ(i))$.
  \item $Am_s(\kk/\QQ(i))\subseteq kerJ_{\k}$.
\end{enumerate}
\end{mainthem}
The  proof of this theorem is based on several results of  units, the class-group  of $\kk$ and
its subgroup of the strongly ambiguous classes.
\section{\textbf{$F.S.U$ OF SOME CM-FIELDS}}
Let us  first collect some results that will be useful in what follows.\par
 Let $m$ and $n$ be two positive square-free integers,  such that $(m,n)=1$; let $\varepsilon_1$ (resp.  $\varepsilon_2$,
$\varepsilon_3$) denote the fundamental unit of $\QQ(\sqrt m)$ (resp.  $\QQ(\sqrt n)$, $\QQ(\sqrt {mn})$).
Put $K_0=\QQ(\sqrt m,\sqrt n)$, $K=K_0(i)$.\\
 \indent Put $B=\{\varepsilon_1, \varepsilon_2, \varepsilon_3, \varepsilon_1\varepsilon_2,
 \varepsilon_1\varepsilon_3, \varepsilon_2\varepsilon_3, \varepsilon_1\varepsilon_2\varepsilon_3\}$
 and\\ $B'=B\cup\{\sqrt \mu/\mu\in B\ and \ \sqrt \mu \in K_0\}$, then  a $F.S.U$ of $K_0$
is a system consisting of three elements chosen from  $B'$ (see \cite{Wa-66} for details).\par
To determine a $F.S.U$ of $K$,  we will use the following result \cite[p.18]{Az-99} that the first author deduced from a theorem of Hasse \cite[\S 21, Satz 15 ]{Ha-52}.
\begin{lem}\label{3}
Let $n\geq2$ be an integer   and $\xi_n$ a $2^n$-th primitive root of unity, then
$$\begin{array}{lllr}\xi_n = \dfrac{1}{2}(\mu_n + \lambda_ni), &\hbox{ where }&\mu_n =\sqrt{2 +\mu_{n-1}},&\lambda_n =\sqrt{2 -\mu_{n-1}}, \\
                               &    & \mu_2=0, \lambda_2=2 &\hbox{ and\quad } \mu_3=\lambda_3=\sqrt 2.
                                              \end{array}
                                            $$
Let $n_0$ be the greatest integer such that  $\xi_{n_0}$ is contained in $K$, $\{\varepsilon'_1, \varepsilon'_2, \varepsilon'_3\}$ a $F.S.U$ of $K_0$ and $\varepsilon$ a unit of  $K_0$ such that ($2 + \mu_{n_0} )\varepsilon$ is a square in  $K_0$ (if it exists). Then a $F.S.U$ of  $K$ is one of the following systems:
\begin{enumerate}[\rm\indent(a)]
\item $\{\varepsilon'_1, \varepsilon'_2, \varepsilon'_3\}$ if $\varepsilon$ does not exist;
\item $\{\varepsilon'_1, \varepsilon'_2, \sqrt{\xi_{n_0}\varepsilon}\}$ if $\varepsilon$
exists;  in this case $\varepsilon = {\varepsilon'_1}^{i_1}{
\varepsilon'_2}^{i_2}\varepsilon'_3$, where $i_1$, $i_2\in \{0, 1\}$ (up to a permutation).
\end{enumerate}
\end{lem}
\begin{lem} \label{2}  If $\varepsilon_1$, $\varepsilon_2$ and $\varepsilon_3$
have negative norms, then
\begin{enumerate}[\upshape\indent(1)]
  \item  If $\varepsilon_1\varepsilon_2\varepsilon_3$ is a square in  $K_0$, then
$\{\varepsilon_1, \varepsilon_2,\sqrt{\varepsilon_1\varepsilon_2\varepsilon_3}\}$ is a $F.S.U$ of $K_0$ and $Q_K=1$.
  \item If $\varepsilon_1\varepsilon_2\varepsilon_3$ is not a square in  $K_0$, then  \{$\varepsilon_1$, $\varepsilon_2$, $\varepsilon_3$\} is a $F.S.U$ of $K_0$ and $Q_K=2$
if and only if  $2\varepsilon_1\varepsilon_2\varepsilon_3$ is a square in $K_0$.
  \item If $Q_K=2$, then $\{\varepsilon_1, \varepsilon_2,\sqrt{i\varepsilon_1\varepsilon_2\varepsilon_3}\}$
  is a $F.S.U$ of $K$.
\end{enumerate}
\end{lem}
\begin{proof}
See Propositions 15, 16 of \cite{Az-05}.
\end{proof}
\begin{lem}[\cite{Az-00}, {Lemma 7}] \label{4}
Let  $p$ be an odd  prime and $\varepsilon_{2p}=x+y\sqrt{2p}$.
If  $N(\varepsilon)=1$, then $x\pm1$ is a square in $\NN$ and $2\varepsilon$
is a square in $\QQ(\sqrt{2p})$.
\end{lem}
\begin{lem}[\cite{Az-00}, {Lemma 5}]\label{5}
Let $d$ be a square-free integer and $\varepsilon_d=x+y\sqrt d$,
where $x$, $y$ are  integers or semi-integers. If $N(\varepsilon)=1$, then $2(x\pm1)$ and
$2d(x\pm1)$ are not squares in  $\QQ$.
\end{lem}
\begin{lem}[\cite{Az-99}, {3.(1) p.19}]\label{6}
Let $d$ be a square-free integer, different from $2$ and $k=\QQ(\sqrt d,i)$, put $\varepsilon_d=x+y\sqrt d$.
\begin{enumerate}[(i)]
  \item If $N(\varepsilon_d)=-1$, then $\{\varepsilon_d\}$ is a $F.S.U$ of $k$.
  \item If $N(\varepsilon_d)=1$, then $\{\sqrt{i\varepsilon_d}\}$ is a $F.S.U$ of $k$ if
   and only if $x\pm1$ is a square in $\NN$ (i.e. $2\varepsilon_d$ is a square in $\QQ(\sqrt d)$). Else $\{\varepsilon_d\}$ is a $F.S.U$ of $k$ (this result is also in \cite{Kub-56}).
\end{enumerate}
\end{lem}
\begin{lem}[\cite{AzTa-08}, {Corollary 3.2}]\label{25}
Let $d$ be square-free integer and   $k=\QQ(\sqrt d,i)$, then $Q_k = 1$ if one of the following  conditions holds:
\begin{enumerate}[\upshape\indent(i)]
  \item $d\equiv1$ $\pmod4$.
  \item There exists an integer  $d'$ dividing $d$ such that $d'\equiv5\pmod 8.$
\end{enumerate}
\end{lem}
\subsection{\textbf{$F.S.U$ of the field $\KK=\QQ(\sqrt{p_1},\sqrt{2p_2},i)$}}
Let $\KK=\kk(\sqrt{p_1})=\\\QQ(\sqrt{p_1}, \sqrt{2p_2}, i)$. We denote by $\varepsilon_1$ (resp. $\varepsilon_2$, $\varepsilon_3$ ) the fundamental unit of  $\QQ(\sqrt{p_1})$ (resp. $\QQ(\sqrt{2p_2})$, $\QQ(\sqrt{2p_1p_2})$) and put  $\varepsilon_3=x+y\sqrt{2p_1p_2}$. The purpose of this sub-paragraph is to establish the following theorem.
 \begin{them}\label{53}
Keep the notations previously mentioned, then
\begin{enumerate}[\upshape\indent(1)]
\item If $N(\varepsilon_2)=N(\varepsilon_3)=-1$, then
\begin{enumerate}[\upshape\indent(i)]
  \item If $\varepsilon_1\varepsilon_2\varepsilon_3$ is a square in $\KK^+$, then $\{\varepsilon_1, \varepsilon_2,\sqrt{\varepsilon_1\varepsilon_2\varepsilon_3}\}$ is a $F.S.U$ of $\KK^+$, $\KK$ and $Q_{\KK}=1$.
 \item Else $\{\varepsilon_1, \varepsilon_2, \varepsilon_3\}$ is a $F.S.U$ of $\KK^+$ and that of $\KK$ is\\ $\{\varepsilon_1, \varepsilon_2,\sqrt{i\varepsilon_1\varepsilon_2\varepsilon_3}\}$  and $Q_{\KK}=2$.
\end{enumerate}
\item If  $N(\varepsilon_2)=-N(\varepsilon_3)=1$, then  the $F.S.U$ of $\KK^+$ is \{$\varepsilon_1$, $\varepsilon_2$, $\varepsilon_3$\} and that of $\KK$ is \{$\varepsilon_1$, $\sqrt{i\varepsilon_2}$, $\varepsilon_3$\} and $Q_{\KK}=2$.
\item If  $N(\varepsilon_3)=-N(\varepsilon_2)=1$, then
\begin{enumerate}[\upshape\indent(i)]
\item If $2p_1(x\pm1)$ is a square in $\NN$, then \{$\varepsilon_1$, $\varepsilon_2$, $\sqrt{\varepsilon_3}$\} is a $F.S.U$ of $\KK^+$,  $\KK$ and $Q_{\KK}=1$.
  \item  Else \{$\varepsilon_1$, $\varepsilon_2$, $\varepsilon_3$\} is a $F.S.U$ of $\KK^+$ and that of $\KK$ is \{$\varepsilon_1$, $\varepsilon_2$, $\sqrt{i\varepsilon_3}$\} and $Q_{\KK}=2$.
\end{enumerate}
\item If  $N(\varepsilon_3)=N(\varepsilon_2)=1$,  then
\begin{enumerate}[\upshape\indent(i)]
\item If $2p_1(x\pm1)$ is a square in $\NN$, then $\{\varepsilon_1, \varepsilon_2, \sqrt{\varepsilon_3}\}$ is a $F.S.U$ of $\KK^+$ and that of $\KK$ is
$\{\varepsilon_1, \sqrt{i\varepsilon_2}, \sqrt{\varepsilon_3}\}$   and $Q_{\KK}=2$.
\item Else $\{\varepsilon_1, \varepsilon_2, \sqrt{\varepsilon_2\varepsilon_3}\}$ is a $F.S.U$ of  $\KK^+$ and that of $\KK$ is
$\{\varepsilon_1, \sqrt{\varepsilon_2\varepsilon_3}, \sqrt{i\varepsilon_3}\}$  and $Q_{\KK}=2$.
\end{enumerate}
\end{enumerate}
\end{them}
\begin{proof}
 See Propositions  \ref{7}, \ref{11}, \ref{12} and \ref{13} below.
  \end{proof}
 \begin{rema}
 Our results in this theorem about unit index of $\KK$ are similar to those  in Theorem 1 (p. 347) of \cite{HiYo}.
 \end{rema}
\begin{propo}\label{7}
Assume that $N(\varepsilon_2)=N(\varepsilon_3)=-1$, then
\begin{enumerate}[\upshape\indent(i)]
  \item If $\varepsilon_1\varepsilon_2\varepsilon_3$ is a square in $\KK^+$, then $\{\varepsilon_1, \varepsilon_2,\sqrt{\varepsilon_1\varepsilon_2\varepsilon_3}\}$ is a $F.S.U$ of $\KK^+$ and $\KK$.
  \item  Else $\{\varepsilon_1, \varepsilon_2, \varepsilon_3\}$ is a $F.S.U$ of $\KK^+$ and $\{\varepsilon_1, \varepsilon_2,\sqrt{i\varepsilon_1\varepsilon_2\varepsilon_3}\}$ is a $F.S.U$ of $\KK$.
\end{enumerate}
\end{propo}
\begin{proof}
(i) If $\varepsilon_1\varepsilon_2\varepsilon_3$ is a square in $\KK^+$, then   Lemma \ref{2} (1) yields the result.\par
(ii) Assume that $\varepsilon_1\varepsilon_2\varepsilon_3$ is not a square in $\KK^+$, then from   Lemma \ref{2} (2) $\{\varepsilon_1, \varepsilon_2, \varepsilon_3\}$ is a $F.S.U$ of $\KK^+$. It remains to determine the $F.S.U$ of $\KK$.\\
As $p_1\equiv p_2\equiv1\pmod4$, then there exist $\pi_1$, $\pi_2$, $\pi_3$ and $\pi_4$ in $\ZZ[i]$ such  that $p_1=\pi_1\pi_2$, $p_2=\pi_3\pi_4$, $\overline{\pi}_1=\pi_2$ and $\overline{\pi}_3=\pi_4$ (the complex conjugate). Let $\varepsilon_1 = a + b \sqrt {p_1}$, where $a$, $b$ are integers or semi-integers.\\
 \indent (a) Suppose that $a$ and $b$ are integers. As  $N(\varepsilon_1)=-1$, then $$(a-i)(a+i)=p_1b^2$$ and the gcd of $a-i$ and $a+i$ divides $2$, so there exist  $b_1$ and $b_2$ in  $\ZZ[i]$ such that $b=b_1b_2$ and
$\left\{\begin{array}{rl}
    a+i &=b_1^2\pi_1,\\
    a-i &=b_2^2\pi_2,
    \end{array}\right.
\text{   or   } \left\{\begin{array}{rl}
    a+i &=ib_1^2\pi_1,\\
    a-i &=-ib_2^2\pi_2;
    \end{array}\right.$\\
 therefore $2a=b_1^2\pi_1+b_2^2\pi_2$ or $2a=ib_1^2\pi_1-ib_2^2\pi_2$, hence
 \begin{center}
 $2\varepsilon_1=(b_1\sqrt{\pi_1}+b_2\sqrt{\pi_2})^2$  or
 $2\varepsilon_1=(b_1\sqrt{i\pi_1}+b_2\sqrt{-i\pi_2})^2$,
 \end{center}
 so
  \begin{center}$\sqrt{2\varepsilon_1}=b_1\sqrt{\pi_1}+b_2\sqrt{\pi_2}$ or
 $2\sqrt{\varepsilon_1}=b_1(1+i)\sqrt{\pi_1}+b_2(1-i)\sqrt{\pi_2}$,
  \end{center}
 we conclude that
   \begin{equation}\label{8}\left\{
   \begin{aligned}\sqrt{2\pi_1\varepsilon_1}& = b_1\pi_1+b_2\sqrt{p_1}\in\KK,\ and  \\
   \sqrt{2\pi_2\varepsilon_1}& =  b_1\sqrt{p_1}+b_2\pi_2\in\KK,
    \ or \\
     2\sqrt{\pi_1\varepsilon_1}& =  b_1(1+i)\pi_1+b_2(1-i)\sqrt{p_1})\in\KK\  and\\
     2\sqrt{\pi_2\varepsilon_1}& = b_1(1+i)\sqrt{p_1}+b_2(1-i)\pi_2\in\KK.
     \end{aligned}
   \right.
\end{equation}
 \indent (b) Let $\varepsilon_1=\frac{1}{2}(a+b\sqrt {p_1})$, where $a$, $b$ are integers, then  \begin{center}$(a-2i)(a+2i)=\pi_1\pi_2b^2$.\end{center} Proceeding as previously we get the same results.\\
 \indent (c) Let $\varepsilon_2=\alpha+\beta\sqrt{2p_2}$, where $\alpha$, $\beta$ are integers, we also find that:
 \begin{equation}\label{9}
 \left\{
 \begin{array}{ll}
\sqrt{2(1+i)\pi_3\varepsilon_2}=\beta_1(1+i)\pi_3+\beta_2\sqrt{2p_2}\in\KK  \text{  and  } \\
\sqrt{2(1-i)\pi_4\varepsilon_2}=\beta_1\sqrt{2p_2}+\beta_2(1-i)\pi_4\in\KK  \text{  or  }\\
\sqrt{(1+i)\pi_3\varepsilon_2}=\frac{1}{2}(\beta_1(1+i)(1\pm i)\pi_3+\beta_2(1\mp i)\sqrt{2p_2})\in\KK
\text{  and  }\\ \sqrt{(1-i)\pi_4\varepsilon_2}=\frac{1}{2}(\beta_1(1\pm i)\sqrt{2p_2}+\beta_2(1-i)(1\mp i)\pi_4)
\in\KK.
 \end{array}\right.
\end{equation}
 \indent (d) Applying the same argument to $\varepsilon_3$, then we  get
  \begin{equation}\label{10}
 \left\{
 \begin{array}{ll}
    \sqrt{2(1+i)\pi_1\pi_3\varepsilon_3}=y_1(1+i)\pi_1\pi_3+y_2\sqrt{2p_1p_2}\in\KK \text{  and  }
    \\ \sqrt{2(1-i)\pi_2\pi_4\varepsilon_3}=y_1\sqrt{2p_1p_2}+y_2(1-i)\pi_2\pi_4 \in\KK \text{ or }\\
    \sqrt{(1+i)\pi_1\pi_3\varepsilon_3}=\\ \frac{1}{2}(y_1(1+i)(1\pm i)\pi_1\pi_3+y_2(1\mp i)\sqrt{2p_1p_2}) \in\KK \
    and
    \\ \sqrt{(1-i)\pi_2\pi_4\varepsilon_3}=\\ \frac{1}{2}(y_1(1\pm i)\sqrt{2p_1p_2}+y_2(1\mp i)(1-i)\pi_2\pi_4)\in\KK.
 \end{array}\right.
\end{equation}
By multiplying the results of equalities (\ref{8}), (\ref{9}) and (\ref{10}), we get
\begin{center}
$\sqrt{\varepsilon_1\varepsilon_2\varepsilon_3}\in \KK^+$ or
$\sqrt{2\varepsilon_1\varepsilon_2\varepsilon_3}\in \KK^+$.
\end{center}
 Finally, note that
$\sqrt{\varepsilon_1\varepsilon_2\varepsilon_3}$ and $\sqrt{2\varepsilon_1\varepsilon_2\varepsilon_3}$
are not both in $\KK^+$, since $\sqrt 2\in \KK^+$.
The rest is a simple deduction from  Lemma \ref{2}(2).
\end{proof}
\begin{propo}\label{11}
If  $N(\varepsilon_2)=-N(\varepsilon_3)=1$, then  the $F.S.U$ of $\KK^+$ is \{$\varepsilon_1$, $\varepsilon_2$, $\varepsilon_3$\} and that of $\KK$ is \{$\varepsilon_1$, $\sqrt{i\varepsilon_2}$, $\varepsilon_3$\}.
\end{propo}
\begin{proof}
\indent As $N(\varepsilon_2)=1$, then, from  Lemma \ref{4}, $a\pm1$ is a square
in $\NN$ and $2\varepsilon_2$ is a square in $\QQ(\sqrt{2p_2})$ ($\varepsilon_2=a+b\sqrt{2p_2}$).
 Note that $\varepsilon_2$ is not a square in $\KK^+$, else we get  $\sqrt 2\in \KK^+$, which is false.\\
\indent  Since $N(\varepsilon_1)=N(\varepsilon_3)=-1$, then $\varepsilon_1$, $\varepsilon_3$ are not squares in $\KK^+$;
  similarly $\varepsilon_1\varepsilon_2$, $\varepsilon_1\varepsilon_3$, $\varepsilon_2\varepsilon_3$
  and $\varepsilon_1\varepsilon_2\varepsilon_3$  are not squares in $\KK^+$,  else we will find that
  $i\in \KK^+$, which is absurd.
Therefore \{$\varepsilon_1$, $\varepsilon_2$, $\varepsilon_3$\} is the $F.S.U$ of $\KK^+$ and as
$2\varepsilon_2$ is a square in $\KK^+$, then, from Lemma \ref{3},  \{$\varepsilon_1$,
$\sqrt{i\varepsilon_2}$, $\varepsilon_3$\} is the $F.S.U$ of $\KK$, which implies that $Q_{\KK}=2$.
\end{proof}
\begin{propo}\label{12}
Assume that  $N(\varepsilon_3)=-N(\varepsilon_2)=1$,   then
\begin{enumerate}[\upshape\indent(i)]
  \item If $2p_1(x\pm1)$ is a square in $\NN$, then \{$\varepsilon_1$, $\varepsilon_2$, $\sqrt{\varepsilon_3}$\} is a $F.S.U$ of $\KK^+$ and $\KK$.
  \item   Else \{$\varepsilon_1$, $\varepsilon_2$, $\varepsilon_3$\} is a $F.S.U$ of $\KK^+$ and that of $\KK$ is \{$\varepsilon_1$, $\varepsilon_2$, $\sqrt{i\varepsilon_3}$\}.
\end{enumerate}
\end{propo}
\begin{proof}
Since $N(\varepsilon_1)=N(\varepsilon_2)=-1$ and $N(\varepsilon_3)=1$, then only $\varepsilon_3$
can be a square in $\KK^+$. The equality $N(\varepsilon_3)=1$ implies  $$(x-1)(x+1)=2p_1p_2y^2,$$ and since  $2p_1p_2\equiv2\pmod4$, thus $2\nmid x$, $2\mid y$ and $(x+1,x-1)=2$, hence with $A=(x\pm1)/2$, $B=(x\mp1)/2$ and $y=2z$ we get $$AB=2p_1p_2z^2,$$
without loss of generality we may assume $2\mid A$. Thus according to Lemma \ref{5} and  to the decomposition uniqueness in $\ZZ$,  we get three possibilities:\\
Case 1: If $p_1\nmid A$, $p_2\nmid A$, then there exist $z_1$, $z_2$ in $\NN$ with $z=z_1z_2$ such that
$$A=2z_1^2, \quad B=p_1p_2z_2^2,$$
hence
$$x\pm1=(2z_1)^2,\quad x\mp1=2p_1p_2z_2^2.$$
Therefore
 $\sqrt{2\varepsilon_3}=\frac{1}{2}(2z_1+z_2\sqrt{2p_1p_2})\in\QQ(\sqrt{2p_1p_2})$ and $\sqrt{\varepsilon_3}\not\in\KK^+$.\\
Case 2: If $p_1\mid A$, $p_2\nmid A$, then  we get
$$A=2p_1z_1^2, \quad B=p_2z_2^2,$$
hence
$$x\pm1=p_1(2z_1)^2,\quad x\mp1=2p_2z_2^2.$$
Therefore
 $\sqrt{2\varepsilon_3}=2z_1\sqrt{p_1}+z_2\sqrt{2p_2}\in\KK^+$ and $\sqrt{\varepsilon_3}\not\in\KK^+$.\\
Case 3:  If $p_1\nmid A$, $p_2\mid A$, then  we get
$$A=2p_2z_1^2, \quad B=p_1z_2^2,$$
hence
$$x\pm1=p_2(2z_1)^2,\quad x\mp1=2p_1z_2^2.$$
Therefore
 $\sqrt{\varepsilon_3}=z_1\sqrt{2p_2}+z_2\sqrt{p_1}\in\KK^+$.\\
 As a result if  $2p_1(x\mp1)$ is a square in $\NN$, then $\{\varepsilon_1, \varepsilon_2, \sqrt{\varepsilon_3}\}$
  is a  $F.S.U$ of $\KK^+$,  $\KK$ and thus $Q_\KK=1$; else  $\sqrt{\varepsilon_3}\not\in\KK^+$  and $\sqrt{2\varepsilon_3}\in\KK^+$,
  this  yields that \{$\varepsilon_1$, $\varepsilon_2$, $\varepsilon_3$\} is a $F.S.U$ of $\KK^+$ and from
  Lemma \ref{3},  \{$\varepsilon_1$, $\varepsilon_2$, $\sqrt{i\varepsilon_3}$\} is a $F.S.U$ of $\KK$, so $Q_\KK=2$.
\end{proof}
\begin{propo}\label{13}
Assume that  $N(\varepsilon_3)=N(\varepsilon_2)=1$, then
\begin{enumerate}[\upshape\indent(i)]
  \item If $2p_1(x\pm1)$ is a square in $\NN$, then $\{\varepsilon_1, \varepsilon_2, \sqrt{\varepsilon_3}\}$ is a $F.S.U$ of $\KK^+$ and that of $\KK$ is $\{\varepsilon_1, \sqrt{i\varepsilon_2}, \sqrt{\varepsilon_3}\}$.
  \item Else $\{\varepsilon_1, \varepsilon_2, \sqrt{\varepsilon_2\varepsilon_3}\}$ is a $F.S.U$ of $\KK^+$ and that of $\KK$ is $\{\varepsilon_1, \sqrt{\varepsilon_2\varepsilon_3}, \sqrt{i\varepsilon_3}\}$.
\end{enumerate}
\end{propo}
\begin{proof}
Since $N(\varepsilon_2)=1$, then, from Lemma \ref{4},  $2\varepsilon_2$ is a square in $k_2$; hence $\varepsilon_2$ is not a square in $\KK^+$.\\
\indent Proceeding as in Proposition \ref{12} we get three cases.\\
\indent (a) If $x\pm1$ is a square in $\NN$, then $\sqrt{\varepsilon_3}\not\in\KK^+$ and $\sqrt{2\varepsilon_3}\in\KK^+$, so $\sqrt{\varepsilon_2\varepsilon_3}\in\KK^+$; hence $\{\varepsilon_1, \varepsilon_2, \sqrt{\varepsilon_2\varepsilon_3}\}$ is a $F.S.U$ of $\KK^+$ and by Lemma \ref{3},  $\{\varepsilon_1, \sqrt{\varepsilon_2\varepsilon_3}, \sqrt{i\varepsilon_3}\}$ is a $F.S.U$ of $\KK$, thus $Q_\KK=2$. It should be noted that, by Lemma \ref{3}, we can take as a $F.S.U$ of $\KK$ one of the following systems $\{\varepsilon_1, \sqrt{\varepsilon_2\varepsilon_3}, \sqrt{i\varepsilon_2}\}$,  $\{\varepsilon_1, \sqrt{i\varepsilon_2}, \sqrt{i\varepsilon_3}\}$.\\
 \indent (b) If  $2p_1(x\pm1)$ is a square in $\NN$, then $\sqrt{\varepsilon_3}\in\KK^+$, so $\sqrt{2\varepsilon_2\varepsilon_3}\in\KK^+$. Thus $\{\varepsilon_1, \varepsilon_2, \sqrt{\varepsilon_3}\}$  is a $F.S.U$ of $\KK^+$ and by Lemma \ref{3}, $\{\varepsilon_1, \sqrt{i\varepsilon_2}, \sqrt{\varepsilon_3}\}$ is a $F.S.U$ of $\KK$, therefore $Q_\KK=2$.\\
 \indent (c) If $p_1(x\pm1)$ is a square in $\NN$, then $\sqrt{\varepsilon_3}\not\in\KK^+$ and $\sqrt{2\varepsilon_3}\in\KK^+$, so $\sqrt{\varepsilon_2\varepsilon_3}\in\KK^+$. The rest is as the case  (a).
\end{proof}
\subsection{\textbf{$F.S.U$ OF THE FIELD $\KK=\QQ(\sqrt2,\sqrt{p_1p_2},i)$}}
Let $\KK=\kk(\sqrt{2})=\\\QQ(\sqrt 2, \sqrt{p_1p_2}, i)$.
Denote by $\varepsilon_1$ (resp. $\varepsilon_2$, $\varepsilon_3$ ) the  fundamental  unit of
$\QQ(\sqrt 2)$ (resp. $\QQ(\sqrt{p_1p_2})$, $\QQ(\sqrt{2p_1p_2})$). Put $I=\{0, 1\}$ and $\varepsilon_3=x+y\sqrt{2p_1p_2}$. Our aim  in this subsection is to state the following theorem,  but  first let us show the lemma.
\begin{lem}\label{14}
 Put $\varepsilon_2=a+b\sqrt{p_1p_2}$; if  $N(\varepsilon_2)=1$, then $a\pm1$ is not a square in $\NN.$
\end{lem}
\begin{proof}
 As $p_1p_2\equiv 1 \pmod 4$, then, from Lemma \ref{25},  the unit index  of $\QQ(\sqrt{p_1p_2}, i)$ is equal to 1; since  $N(\varepsilon_2)=1$, so  assertion 3.(1).(ii) on p.19 of \cite{Az-99},  yields that  $a\pm1$ is not a square in $\NN$.
\end{proof}
\begin{them}\label{54}
Keep notations  mentioned above,  Then $Q_{\KK}=1$ and
\begin{enumerate}[\upshape\indent(1)]
\item If $N(\varepsilon_2)=N(\varepsilon_3)=-1$, then
\begin{enumerate}[\upshape\indent(i)]
  \item If $\varepsilon_1\varepsilon_2\varepsilon_3$ is a square in $\KK^+$, then $\{\varepsilon_1, \varepsilon_2,\sqrt{\varepsilon_1\varepsilon_2\varepsilon_3}\}$ is a $F.S.U$ of $\KK^+$, $\KK$.
 \item Else $\{\varepsilon_1, \varepsilon_2, \varepsilon_3\}$ is a $F.S.U$ of $\KK^+$, $\KK$.
\end{enumerate}
\item If  $N(\varepsilon_2)=-N(\varepsilon_3)=1$, then  the $F.S.U$ of $\KK^+$,  $\KK$ is \{$\varepsilon_1$, $\varepsilon_2$, $\varepsilon_3$\}.
\item If  $N(\varepsilon_3)=-N(\varepsilon_2)=1$, then
\begin{enumerate}[\upshape\indent(i)]
\item If $x\pm1$ is a square in $\NN$, then \{$\varepsilon_1$, $\varepsilon_2$, $\sqrt{\varepsilon_3}$\} is a $F.S.U$ of $\KK^+$,  $\KK$.
  \item  Else \{$\varepsilon_1$, $\varepsilon_2$, $\varepsilon_3$\} is a $F.S.U$ of $\KK^+$, $\KK$.
\end{enumerate}
\item If  $N(\varepsilon_3)=N(\varepsilon_2)=1$,  then
\begin{enumerate}[\upshape\indent(i)]
\item If $x\pm1$ is a square in $\NN$, then $\{\varepsilon_1, \varepsilon_2, \sqrt{\varepsilon_3}\}$ is a $F.S.U$ of $\KK^+$, $\KK$.
\item Else $\{\varepsilon_1, \varepsilon_2, \sqrt{\varepsilon_2\varepsilon_3}\}$ is a $F.S.U$ of $\KK^+$, $\KK$.
\end{enumerate}
\end{enumerate}
\end{them}
\begin{proof}
See  Propositions \ref{15}, \ref{16}, \ref{21} and \ref{22} below.
\end{proof}
\begin{rema}
The unit index of $\KK$ is always equal to 1, which is compatible with theorem 2 of \cite{HiYo}.
\end{rema}
\begin{propo}\label{15}
Assume  that $N(\varepsilon_2)=-N(\varepsilon_3)=1$,
then $\{\varepsilon_1,  \varepsilon_2,  \varepsilon_3 \}$ is a $F.S.U$ of  $\KK^+$ and  $\KK$.
\end{propo}
\begin{proof}
Since $\varepsilon_1$, $\varepsilon_3$ have   negative norms, then they are not squares in $\KK^+$,
similarly $\varepsilon_1\varepsilon_2$, $\varepsilon_1\varepsilon_3$, $\varepsilon_2\varepsilon_3$ and
$\varepsilon_1\varepsilon_2\varepsilon_3$  are not squares in $\KK^+$, else  by taking a suitable norm we
get $i\in\KK^+$, which is false. Furthermore
 $(2+\sqrt2)\varepsilon_1^i\varepsilon_2^j\varepsilon_3^k$ cannot be a square in $\KK^+$, for all $i$, $j$ and $k$ of
 I, as otherwise with some $\alpha\in\KK^+$ we would have $\alpha^2=(2+\sqrt 2)\varepsilon_1^i\varepsilon_2^j\varepsilon_3^k$, so $N_{\KK^+/k_2}^2(\alpha)=2(-1)^{i+k}\varepsilon_2^{2j}$,  yielding that $\sqrt2\in\QQ(\sqrt{p_1p_2})$, which is absurd.\\
 \indent  If $\varepsilon_2=a+b\sqrt{p_1p_2}$, then $a^2-1=b^2p_1p_2$. Proceeding as in Proposition \ref{12} and taking into account   Lemma \ref{5},  we get the following cases.\\
 \indent (i) If $p_1p_2(a\pm1)$ is a square in $\NN$, then there exists $(b_1, b_2)\in\ZZ^2$ such that
\begin{center} $a\mp1=b_1^2 \text{ and }
 a\pm1=b_2^2p_1p_2, $\end{center}
 so $a\mp1$ is a square in $\NN$, which contradicts  Lemma \ref{14}.\\
\indent (ii) If $p_1(a\pm1)$ is a square in $\NN$,
then there exists  $(b_1, b_2)\in\ZZ^2$ such that
\begin{center}
 $a\pm1=p_1b_1^2 \text{ and }
 a\mp1=p_2b_2^2,$\end{center}
thus $\sqrt{\varepsilon_2}=\frac{1}{2}(b_1\sqrt{2p_1}+b_2\sqrt{2p_2})$, so $\sqrt{\varepsilon_2}\not\in\KK^+$, $\sqrt{p_1\varepsilon_2}\in\KK^+$ and $\sqrt{p_2\varepsilon_2}\in\KK^+$.\\
\indent (iii)  If $2p_1(a\pm1)$ is a square in $\NN$, then the same argument shows\\ that $\sqrt{\varepsilon_2}=b_1\sqrt{p_1}+b_2\sqrt{p_2}$, $\sqrt{p_1\varepsilon_2}\in\KK^+$ and $\sqrt{p_2\varepsilon_2}\in\KK^+$.\par
Therefore we deduce that $\{\varepsilon_1,  \varepsilon_2,  \varepsilon_3 \}$ is a $F.S.U$ of $\KK^+$ and from Lemma \ref{3}, $\KK$ has the same $F.S.U$.
\end{proof}
\begin{propo}\label{16}
Assume that  $N(\varepsilon_2)=N(\varepsilon_3)=-1$, then
\begin{enumerate}[\upshape\indent(i)]
  \item If $\varepsilon_1\varepsilon_2\varepsilon_3$ is a square in $\KK_3^+$, then
 $\{ \varepsilon_1, \varepsilon_2, \sqrt{\varepsilon_1\varepsilon_2\varepsilon_3}\}$ is a $F.S.U$ of $\KK_3^+$ and $\KK_3$.
  \item  Else $\{ \varepsilon_1, \varepsilon_2, \varepsilon_3\}$ is a $F.S.U$ of  $\KK_3^+$ and $\KK_3$.
\end{enumerate}
\end{propo}
\begin{proof}
We proceed as in Proposition \ref{15} to prove that  $(2+\sqrt2)\varepsilon_1^i\varepsilon_2^j\varepsilon_3^k$ is not a square in $\KK^+$, for all $i$, $j$ and $k$ of $I$, and we apply  Lemmas \ref{2}, \ref{3} and Remark \ref{17} bellow.
\end{proof}
\begin{rema}\label{17}
 We proceed as in the proof of Proposition \ref{7} to prove the following:\\
 \indent (a) If $\varepsilon_2=x+y\sqrt{p_1p_2}$, then
  \begin{equation}\label{18}\left\{
   \begin{aligned}\sqrt{\varepsilon_2}& = y_1\sqrt{\pi_1\pi_3}+y_2\sqrt{\pi_2\pi_4},\ or  \\
   \sqrt{\varepsilon_2}& = y_1\sqrt{\pi_1\pi_4}+y_2\sqrt{\pi_2\pi_3},\ or
    \  \\
    \sqrt{2\varepsilon_2}& = y_1\sqrt{\pi_1\pi_3}+y_2\sqrt{\pi_2\pi_4},\ or \\
     \sqrt{2\varepsilon_2}& = y_1\sqrt{\pi_1\pi_4}+y_2\sqrt{\pi_2\pi_3},
     \end{aligned}
   \right.
\end{equation}
where  $y_i$ are in $\ZZ[i]$ or $\frac{1}{2}\ZZ[i]$.\\
  \indent (b) If  $\varepsilon_3=a+b\sqrt{2p_1p_2}$, then
    \begin{equation}\label{19}
    \left\{
   \begin{aligned}
  \sqrt{\varepsilon_3}& = b_1\sqrt{(1+i)\pi_1\pi_3}+b_2\sqrt{(1-i)\pi_2\pi_4}),\ or  \\
  \sqrt{\varepsilon_3}& = b_1\sqrt{(1+i)\pi_1\pi_4}+b_2\sqrt{(1-i)\pi_2\pi_3}),\ or
    \  \\
    \sqrt{2\varepsilon_3}& = b_1\sqrt{(1+i)\pi_1\pi_3}+b_2\sqrt{(1-i)\pi_2\pi_4}),\ or \\
    \sqrt{2\varepsilon_3}& = b_1\sqrt{(1+i)\pi_1\pi_3}+b_2\sqrt{(1-i)\pi_2\pi_4}),
     \end{aligned}
   \right.
\end{equation}
where $b_i$ are in  $\ZZ[i]$ or $\frac{1}{2}\ZZ[i]$.\\
  Note  at the end that:
 \begin{equation}\label{20}
  \sqrt{2\varepsilon_1}=\sqrt{1+i}+\sqrt{1-i}.\end{equation}
So by multiplying  results of equalities  $(\ref{18})$, $(\ref{19})$ and $(\ref{20})$ we get
\begin{center}
$\sqrt{\varepsilon_1\varepsilon_2\varepsilon_3}=\alpha+\beta\sqrt{2}+\gamma\sqrt{p_1p_2}+\delta\sqrt{2p_1p_2}
\in\QQ(\sqrt2, \sqrt{p_1p_2})$ or\\   $\sqrt{\varepsilon_1\varepsilon_2\varepsilon_3}=\alpha\sqrt{p_1}+
\beta\sqrt{p_2}+\gamma\sqrt{2p_1}+\delta\sqrt{2p_2}\not\in\QQ(\sqrt2, \sqrt{p_1p_2})$,\end{center} where
$\alpha$, $\beta$, $\gamma$ and $\delta$ are in $\QQ$.
\end{rema}
\begin{propo}\label{21}
Assume that  $N(\varepsilon_3)=-N(\varepsilon_2)=1$,  then
\begin{enumerate}[\upshape\indent(i)]
\item If $x\pm1$ is a square in $\NN$, then $\{\varepsilon_1, \varepsilon_2, \sqrt{\varepsilon_3}\}$ is a $F.S.U$ of $\KK$ and $\KK^+$.
 \item Else $\{ \varepsilon_1, \varepsilon_2, \varepsilon_3\}$ is a $F.S.U$ of $\KK^+$ and $\KK$.
 \end{enumerate}
\end{propo}
\begin{proof}
As the norms of $\varepsilon_1$, $\varepsilon_2$ are negative, then  proceeding as in  Proposition
\ref{15}, we prove that only $\varepsilon_3$ can be a square in $\KK^+$.\\
 \indent (i) According to Lemma \ref{6},  $2\varepsilon_3$ is a square in $\QQ(\sqrt{2p_1p_2})$, and since $\sqrt 2\in\KK^+$, so $\sqrt{\varepsilon_3}\in\KK^+$, which yields that
 $\{\varepsilon_1, \varepsilon_2, \sqrt{\varepsilon_3}\}$ is a $F.S.U$ of $\KK$ and, by Lemma \ref{3}, is also a $F.S.U$ of $\KK^+$.\\
 \indent (ii)  $\varepsilon_3$ is not a square in $\KK^+$, then  $\{ \varepsilon_1, \varepsilon_2, \varepsilon_3\}$
 is a $F.S.U$ of $\KK$, $\KK^+$ (Lemma \ref{3}).
\end{proof}
\begin{propo}\label{22}
Assume that $N(\varepsilon_3)=N(\varepsilon_2)=1$,  then
\begin{enumerate}[\upshape\indent(i)]
 \item If $x\pm1$ a square in $\NN$, then $\{\varepsilon_1, \varepsilon_2, \sqrt{\varepsilon_3}\}$ is a $F.S.U$ of $\KK$ and $\KK^+$.
 \item  Else $\{ \varepsilon_1, \varepsilon_2, \sqrt{\varepsilon_2\varepsilon_3}\}$ is a $F.S.U$ of $\KK^+$ and  $\KK$.
 \end{enumerate}
\end{propo}
\begin{proof}
Since   $N(\varepsilon_1)=-1$, then $\varepsilon_1$, $\varepsilon_1\varepsilon_2$,
$\varepsilon_1\varepsilon_3$,  $\varepsilon_1\varepsilon_2\varepsilon_3$ and
$(2+\sqrt2)\varepsilon_1^i\varepsilon_2^j\varepsilon_3^k$ are not squares in $\KK^+$, for all $i$, $j$
and $k$ in I.\\
\indent As   $N(\varepsilon_2)=1$, then from  Lemma \ref{14}, there exist $y_1$,  $y_2$ in $\ZZ$
such that: \begin{equation}\label{23}\sqrt{2\varepsilon_2}=y_1\sqrt{p_1}+y_2\sqrt{p_2}.\end{equation}
\indent The equality  $N(\varepsilon_3)=1$ and the Lemma \ref{5} imply the existence of $a_1$, $a_2$, $n$ and $m$ in
$\NN$, such that
\begin{center}
$n=p_1 \text{ and }
 m=2p_2 \text{ or }
 n=2p_1 \text{ and }
 m=p_2,$
\text{ and }\\
$\left\{\begin{array}{ll}
 x\pm1=a_1^2,\\
 x\mp1=2p_1p_2,
 \end{array}\right.
  \text{ or }
 \left\{\begin{array}{ll}
 x\pm1=na_1^2,\\
 x\mp1=ma_2^2;
  \end{array}\right.$
 \end{center}
 we get
 \begin{equation}\label{24} \sqrt{\varepsilon_3}=\frac{1}{2}(a_1\sqrt{2}+2a_2\sqrt{p_1p_2})\text{ or }
 \sqrt{\varepsilon_3}=\frac{1}{2}(2a_1\sqrt{n}+a_2\sqrt{m}).\end{equation}
From equalities (\ref{23}) and (\ref{24}) we deduce that if  $x\pm1$ is a square in $\NN$, then
$\sqrt{\varepsilon_3}\in\KK^+$, else  $\sqrt{\varepsilon_2\varepsilon_3}\in\KK^+$. Hence the results.
\end{proof}
\section{\textbf{Some results}}
In what remains of this paper, we adopt the following notations. Let
$p_1=e^2+4f^2=\pi_1\pi_2$, $p_2=g^2+4h^2=\pi_3\pi_4$, $\pi_1=e+2if$, $\pi_2=e-2if$, $\pi_3=g+2ih$, $\pi_4=g-2ih$. Let $\mathcal{H}_j$
be the prime ideal of $\kk$ above $\pi_j$,  hence $\mathcal{H}_j^2=(\pi_j)$, for all $j\in\{1, 2, 3, 4\}$. As $2$ is totally ramified in $\kk$, let $\mathcal{H}_0$ be the prime ideal of $\kk$ above $1+i$, so $\mathcal{H}_0^2=(1+i)$
\begin{propo}\label{26}
Let $d$ be a square-free integer, $k=\QQ(\sqrt d,i)$, $a+ib$ an element of $\ZZ[i]$ and $\mathcal{H}$ an ideal of $k$ such that $\mathcal{H}^2=(a+ib)$. Put $\varepsilon_d=x+y\sqrt d$ the  fundamental unit of $\QQ(\sqrt d)$, so\\
 $(1)$ If $\sqrt{a^2+b^2}\not\in\QQ(\sqrt d)$, then  $\mathcal{H}$ is not principal in $k$.\\
 $(2)$ If $a^2+b^2=d$, then
 \begin{enumerate}[\upshape\indent(a)]
   \item If  $N(\varepsilon_d)=1$, then $\mathcal{H}$ is not  principal in $k$.
   \item If $N(\varepsilon_d)=-1$, then:\\
   (i) If $(ax\pm yd)\pm b$ or  $2(-xb\pm yd)\pm a$ is a square in $\NN$, then $\mathcal{H}$ is  principal in $k$.\\
   (ii) Else $\mathcal{H}$ is not principal in $k$.
 \end{enumerate}
\end{propo}
\begin{proof}
See Proposition $1$ of \cite{AZT12-2}
\end{proof}
\begin{propo}\label{28}
Let $d$ be a composite integer, even, square-free and product at least of three primes. Let $k=\QQ(\sqrt d,i)$, $p$ an odd prime  and $\mathcal{H}$
an ideal of $k$ such that $\mathcal{H}^2=(p)$. Let $\varepsilon_d=x+y\sqrt d$ denote the  fundamental unit of $\QQ(\sqrt d)$. Then
\begin{enumerate}[\rm\indent(1)]
  \item If $N(\varepsilon_d)=-1$, then $\mathcal{H}$ is not  principal in $k$.
  \item  If  $N(\varepsilon_d)=1$, then
\begin{enumerate}[\rm\indent(i)]
  \item  If $\{\varepsilon_d\}$ is $F.S.U$ of $k$, then $\mathcal{H}$ is principal in $k$ if and only if
  $2p(x\pm1)$ or $p(x\pm1)$  is a square in $\NN$.
  \item  Else  $\mathcal{H}$ is not principal in $k$.
\end{enumerate}
\end{enumerate}
\end{propo}
\begin{proof}
See Proposition $2$ of \cite{AZT12-2}
\end{proof}
We finish this paragraph by the  following two Lemmas.
\begin{lem}\label{41}
Let  $d=2p_1p_2$  and $\varepsilon_d=x+y\sqrt d$. If  $N(\varepsilon_d)=-1$, then there exist $y_1$ and $y_2$ in $\ZZ[i]$ such that  $\sqrt{\varepsilon_d}$
  takes one of following forms:
\begin{enumerate}[\rm\indent (1)]
  \item $\frac{1}{2}[y_1(1+i)\sqrt{(1\pm i)\pi_1\pi_3}+y_2(1-i)\sqrt{(1\mp i)\pi_2\pi_4}]$.
 \item $\frac{1}{2}[y_1(1+i)\sqrt{(1\pm i)\pi_1\pi_4}+y_2(1-i)\sqrt{(1\mp i)\pi_2\pi_3}]$.
\end{enumerate}
\end{lem}
\begin{proof}
As  $N(\varepsilon_d)=-1$, then $x^2+1=y^2d$, hence there exist $y_1$, $y_2$ in  $\ZZ[i]$ such that
\begin{center}
$\left\{\begin{array}{ll}
 x\pm i=(1+i)\pi_1\pi_3y_1^2,\\
 x\mp i=(1-i)\pi_2\pi_4y_2^2,
 \end{array}\right.$ or
 $\left\{\begin{array}{ll}
 x\pm i=(1+i)\pi_1\pi_4y_1^2,\\
 x\mp i=(1-i)\pi_2\pi_3y_2^2,
 \end{array}\right.$ or \\
 $\left\{\begin{array}{ll}
 x\pm i=i(1+i)\pi_1\pi_3y_1^2,\\
 x\mp i=-i(1-i)\pi_2\pi_4y_2^2,
 \end{array}\right.$ or
 $\left\{\begin{array}{ll}
 x\pm i=i(1+i)\pi_1\pi_4y_1^2,\\
 x\mp i=-i(1-i)\pi_2\pi_3y_2^2;
 \end{array}\right.$
 \end{center}
 therefore\\ $2x=(1+i)\pi_1\pi_3y_1^2+(1-i)\pi_2\pi_4y_2^2$ or $2x=(1+i)\pi_1\pi_4y_1^2+(1-i)\pi_2\pi_3y_2^2$ or
  $2x=i(1+i)\pi_1\pi_3y_1^2-i(1-i)\pi_2\pi_4y_2^2$ or $2x=i(1+i)\pi_1\pi_4y_1^2-i(1-i)\pi_2\pi_3y_2^2$,\\ so\\
 $\varepsilon_d=\frac{1}{2}[y_1\sqrt{(1+i)\pi_1\pi_3}+y_2\sqrt{(1-i)\pi_2\pi_4}]^2$ or \\
  $\varepsilon_d=\frac{1}{2}[y_1\sqrt{(1+i)\pi_1\pi_4}+y_2\sqrt{(1-i)\pi_2\pi_3}]^2$ or \\
  $\varepsilon_d=\frac{1}{2}[y_1(1+i)\sqrt{\frac{(1+i)\pi_1\pi_3}{2}}+y_2(1-i)\sqrt{\frac{(1-i)\pi_2\pi_4}{2}}]^2$ or\\
   $\varepsilon_d=\frac{1}{2}[y_1(1+i)\sqrt{\frac{(1+i)\pi_1\pi_4}{2}}+y_2(1-i)\sqrt{\frac{(1-i)\pi_2\pi_3}{2}}]^2$,\\ hence\\
 $\sqrt{\varepsilon_d}=\frac{1}{2}[y_1(1+i)\sqrt{(1\pm i)\pi_1\pi_3}+y_2(1-i)\sqrt{(1\mp i)\pi_2\pi_4}]$ or \\
 $\sqrt{\varepsilon_d}=\frac{1}{2}[y_1(1+i)\sqrt{(1\pm i)\pi_1\pi_4}+y_2(1-i)\sqrt{(1\mp i)\pi_2\pi_3}]$.
\end{proof}
\begin{lem}\label{42}
Let   $d=2p_1p_2$  and $\kk=\QQ(\sqrt{d},i)$. If $N(\varepsilon_d)=-1$, then
\begin{enumerate}[\rm\indent(i)]
  \item If  $\sqrt{\varepsilon_d}$ takes the form $(1)$ of Lemma $\ref{41}$, then
   $\mathcal{H}_0\mathcal{H}_1\mathcal{H}_3$ and $\mathcal{H}_0\mathcal{H}_2\mathcal{H}_4$ are principal in $\kk$.
  \item If  $\sqrt{\varepsilon_d}$ takes the form $(2)$ of Lemma $\ref{41}$, then
   $\mathcal{H}_0\mathcal{H}_1\mathcal{H}_4$ and $\mathcal{H}_0\mathcal{H}_2\mathcal{H}_3$ are principal in $\kk$.
\end{enumerate}
\end{lem}
\begin{proof}
It is easy to see that for all $j\in\{1, 2, 3, 4\}$, $\pi_j$ is ramified in $\kk/\QQ(i)$, thus  $\mathcal{H}_j^2=(\pi_j)$. On the other hand, $2$ is totally ramified in $\kk$ and $\mathcal{H}_0^2=(1+i)$. Hence for all $j\in\{0, 1, 2, 3, 4\}$, Proposition \ref{26}, states that  $\mathcal{H}_j$ is not principal in $\kk$, since $\sqrt{p_1}$, $\sqrt{p_2}$ and $\sqrt2$ are not in $\QQ(\sqrt{2p_1p_2})$.\par
 The ideal $\mathcal{H}_0\mathcal{H}_1\mathcal{H}_3$ is principal in $\kk$, if and only if there exists a unit $\varepsilon\in\kk$ such that   \begin{equation}\label{55}(1+i)\pi_1\pi_3\varepsilon=\alpha^2,\end{equation} where $\alpha\in\kk$. As  $N(\varepsilon_d)=-1$,
so, by Lemma \ref{6},  $Q_\kk=1$; hence $\varepsilon$ is  either real or purely imaginary.\par
Put $\alpha=\alpha_1+i\alpha_2$,  with $\alpha_1$,
 $\alpha_2\in\QQ(\sqrt{2p_1p_2})$, and suppose  $\varepsilon$ is  real (same proof if it is purely imaginary), since $\pi_1\pi_3=(e+2if)(g+2ih)=(eg-4fh)+2i(eh+gf)$, then the equation (\ref{55}) is equivalent to
$$\alpha_1^2-\alpha_2^2+2i\alpha_1\alpha_2=\varepsilon[(eg-4fh)-2(eh+fg)]+i\varepsilon_d[(eg-4fh)+2(eh+gf)],$$
hence
$$\begin{array}{ll}
   \alpha_1^2-\alpha_2^2 &=\varepsilon[(eg-4fh)-2(eh+fg)], and \\
   2\alpha_1\alpha_2 & =\varepsilon[(eg-4fh)+2(eh+gf)],
\end{array}$$  so we get
 $$\alpha_2=\frac{\varepsilon[(eg-4fh)+2(eh+gf)]}{2\alpha_1},$$ thus
 $$4\alpha_1^4-4\varepsilon[(eg-4fh)-2(eh+fg)]\alpha_1^2-[(eg-4fh)+2(eh+fg)]^2\varepsilon^2=0,$$
the discriminant of this equation is $\Delta'=4\varepsilon^2d$,
 which implies that $$\alpha_1^2=\frac{\varepsilon}{4}[2[(eg-4fh)-2(eh+fg)]\pm2\sqrt d].$$
 On the other hand, \begin{center} $(1+i)\pi_1\pi_3+(1-i)\pi_2\pi_4 =2(eg-4fh)-4(eh+fg)$ and \\
          $\sqrt d  =\sqrt{(1-i)\pi_1\pi_3}\sqrt{(1+i)\pi_2\pi_4},$ \end{center}
   then
   $$\begin{array}{ll}
   \alpha_1^2&=\frac{\varepsilon}{4}(\sqrt{(1-i)\pi_1\pi_3}+\sqrt{(1+i)\pi_2\pi_4})^2,\text{ so }\\
    \alpha_1&=\frac{\sqrt{\varepsilon}}{2}(\sqrt{(1-i)\pi_1\pi_3}+\sqrt{(1+i)\pi_2\pi_4}),
    \end{array}$$
   therefore  if $\varepsilon=\varepsilon_d$ and $\sqrt{\varepsilon_d}$ takes the value  (1)  of Lemma \ref{41}, we get
 $$
   \alpha_1 =\frac{1}{4}(2y_1\pi_1\pi_3+2y_2\pi_2\pi_4+(y_1(1+i)+y_2(1-i))\sqrt d),$$ and
   $$\alpha_2 =\frac{\varepsilon_d[(eg-4fh)+2(eh+gf)]}{2\alpha_1},
     $$
    and it is easy to see that $\alpha_1$, $\alpha_2$
   $\in\QQ(\sqrt{2p_1p_2})$; hence $\mathcal{H}_0\mathcal{H}_1\mathcal{H}_3$ is principal in $\kk$.
   By writing $\pi_2\pi_4$ in terms of e, f, g and h we prove similarly that $\mathcal{H}_0\mathcal{H}_2\mathcal{H}_4$ is principal.\\
   \indent Proceeding  similarly  we prove that $\mathcal{H}_0\mathcal{H}_1\mathcal{H}_4$ and
  $\mathcal{H}_0\mathcal{H}_2\mathcal{H}_3$ are principal in $\kk$, if
   $\sqrt{\varepsilon_d}$ takes the value  (2)  of Lemma \ref{41}.
\end{proof}
\section{\textbf{The strongly ambiguous classes of $\kk/\QQ(i)$}}
 Let $F=\QQ(i)$ and  Galois$(\kk/F)=\langle\sigma\rangle$. We denote by  $Am(\kk/F)$  the group of the ambiguous classes of $\kk/F$, that are classes of $\kk$ fixed under $\sigma$, we denote also by $Am_s(\kk/F)$ the subgroup of $Am(\kk/F)$ generated by  the strongly ambiguous classes, which are classes  of $\kk$ containing at least one ideal invariant under  $\sigma$. The genus number, $[(\kk/F)^*:\kk]$, is given by the ambiguous class number formula (see \cite{Ch-33}):
 \begin{equation}\label{51}
|Am(\kk/F)|=[(\kk/F)^*:\kk]=\frac{h(F)2^{t-1}}{[E_F: E_F\cap N_{\kk/F}(\kk^\times)]},
\end{equation}
where $h(F)$ is the class number of $F$, $t$ is the number of finite and infinite primes of $F$ ramified in $\kk/F$. Moreover as the class number
of  $F$ is equal to $1$, so it is well known that
 \begin{equation}\label{51}|Am(\kk/F)|=[(\kk/F)^*:\kk]=2^r,\end{equation}
 where $r=\text{rank}\mathbf{C}_{\mathds{k},2}$.
  The relation between  $|Am(\kk/F)|$ and $|Am_s(\kk/F)|$ is given by the formula:
\begin{equation}\label{50}
\frac{|Am(\kk/F)|}{|Am_s(\kk/F)|}=[E_F\cap N_{\kk/F}(\kk^\times):N_{\kk/F}(E_\kk)].
\end{equation}
Since $\mathcal{H}_0^2=(1+i)$ and for all $j\in\{1, 2, 3, 4\}$, $\mathcal{H}_j^2=(\pi_j)$, so for all $j\in\{0, 1, 2, 3, 4\}$, $[\mathcal{H}_j]$ is a strongly ambiguous class of $\kk/F$ i.e. $[\mathcal{H}_j]\in Am_s(\kk/F)$.
\begin{propo}\label{52}
 Let  $d=2p_1p_2$  and $\kk=\QQ(\sqrt{d},i)$.
\begin{enumerate}[\upshape\indent(1)]
\item $\k\varsubsetneq (\kk/F)^*$.
  \item Assume that $(p_1\equiv p_2\equiv1 \text{ mod $8$ }$ and $Q_\kk=1)$ or $(p_1\equiv5$ or $p_2\equiv5 \text{ mod $8$ })$, then
  \begin{enumerate}[\upshape\indent(i)]
    \item If $N(\varepsilon_d)=-1$, then $Am_s(\kk/\QQ(i))=\langle [\mathcal{H}_0], [\mathcal{H}_1], [\mathcal{H}_2]\rangle=\\ \langle [\mathcal{H}_0], [\mathcal{H}_3], [\mathcal{H}_4]\rangle$.
   \item  Else  $Am_s(\kk/\QQ(i))=\langle [\mathcal{H}_0], [\mathcal{H}_1], [\mathcal{H}_3]\rangle$.
  \end{enumerate}
  \item If $p_1\equiv p_2\equiv1 \text{ mod $8$ }$ and $Q_\kk=2$, then \\ $Am_s(\kk/\QQ(i))=\langle [\mathcal{H}_0], [\mathcal{H}_1], [\mathcal{H}_2], [\mathcal{H}_3]\rangle$.
\end{enumerate}
\end{propo}
\begin{proof}
(1) As $\kk=\QQ(\sqrt{2p_1p_2}, i)$, so $[\k:\kk]=4$; moreover, according to \cite[Proposition 2, p.90]{McPaRa-95}, $r=\text{rank}\mathbf{C}_{\mathds{k},2}=4$ if $p_1\equiv p_2\equiv1 \pmod8$ and $r=\text{rank}\mathbf{C}_{\mathds{k},2}=3$ if $p_1\equiv5$ or $p_2\equiv5 \pmod8$, so $[(\kk/F)^*:\kk]=8 \text{ or } 16$, which implies  the result.\\
 (2)  Note first
that if $p_1\equiv5$ or $p_2\equiv5 \pmod8$, then  Lemma \ref{25}  states that $Q_\kk=1$, hence
$N(\varepsilon_d)=-1$ or ($N(\varepsilon_d)=1$ and $x\pm1$ is not a square in $\NN$) (see Lemma \ref{6}),
where $\varepsilon_d=x+y\sqrt{2p_1p_2}$.\par
 (i) Since $(\mathcal{H}_1\mathcal{H}_2)^2=(p_1)$,  $(\mathcal{H}_3\mathcal{H}_4)^2=(p_2)$, so Proposition \ref{28}
implies that  $\mathcal{H}_1\mathcal{H}_2$,  $\mathcal{H}_3\mathcal{H}_4$ are not principal in $\kk$, hence $\mathcal{H}_1$ and $\mathcal{H}_2$ (resp. $\mathcal{H}_3$ and $\mathcal{H}_4$) lie in different  classes. Moreover  Lemma \ref{42} states that
 $[\mathcal{H}_0\mathcal{H}_1]=[\mathcal{H}_3]$
 and  $[\mathcal{H}_0\mathcal{H}_2]=[\mathcal{H}_4]$ or $[\mathcal{H}_0\mathcal{H}_2]=[\mathcal{H}_3]$
and $[\mathcal{H}_0\mathcal{H}_1]=[\mathcal{H}_4]$.
 Since $(\mathcal{H}_0\mathcal{H}_1\mathcal{H}_2)^2=((1+i)p_1)$ (resp. $(\mathcal{H}_0\mathcal{H}_3\mathcal{H}_4)^2=((1+i)p_2)$)
and $\sqrt{2}\not\in\QQ(\sqrt{2p_1p_2})$, then  Proposition \ref{26} yields that
$\mathcal{H}_0\mathcal{H}_1\mathcal{H}_2$  and  $\mathcal{H}_0\mathcal{H}_3\mathcal{H}_4$ are not principal in $\kk$. Therefore
$$\langle [\mathcal{H}_0], [\mathcal{H}_1], [\mathcal{H}_2]\rangle\subseteq Am_s(\kk/\QQ(i)) \text{ and }\langle [\mathcal{H}_0], [\mathcal{H}_3], [\mathcal{H}_4]\rangle\subseteq Am_s(\kk/\QQ(i)).$$
 On the other hand, since $Q_\kk=1$, then by Lemma \ref{6} we get  $E_\kk=\langle i, \varepsilon_3\rangle$, so $N_{\kk/F}(E_\kk)=\langle -1\rangle$.\par
 (a) Suppose $p_1\equiv p_2\equiv1 \pmod8$, so
$r=\text{rank}\mathbf{C}_{\mathds{k},2}=4$ and it is well known that $i$ is norm in $\kk/F$, hence  $|Am(\kk/F)|=2^4$
and $E_F\cap N_{\kk/F}(\kk^\times)=E_F=\langle i\rangle$, therefore formula (\ref{50}) yields that $|Am_s(\kk/F)|=8$,
 this states that  $$Am_s(\kk/\QQ(i))=\langle [\mathcal{H}_0], [\mathcal{H}_1], [\mathcal{H}_2]\rangle=\langle [\mathcal{H}_0], [\mathcal{H}_3], [\mathcal{H}_4]\rangle.$$\par
 (b) Suppose $p_1\equiv5$ or $p_2\equiv5 \pmod8$, so
$r=\text{rank}\mathbf{C}_{\mathds{k},2}=3$ and in this case $i$ is not norm in $\kk/F$, therefore
$|Am(\kk/F)|=2^3$ and $E_F\cap N_{\kk/F}(\kk^\times)=\langle -1\rangle$, so formula (\ref{50}) yields that
$|Am(\kk/F)|=|Am_s(\kk/F)|=8$, and the result derived.\par
(ii) If $N(\varepsilon_3)=1$ and $x\pm1$ is not a square in $\NN$, then from Lemma \ref{5} we get  $p_1(x\pm1)$
and $2p_2(x\mp1)$ or $p_2(x\pm1)$ and $2p_1(x\mp1)$ are squares in $\NN$, hence Proposition \ref{28} implies
that $\mathcal{H}_1\mathcal{H}_2$ and $\mathcal{H}_3\mathcal{H}_4$ are principal in $\kk$, since
$(\mathcal{H}_1\mathcal{H}_2)^2=(p_1)$ and $(\mathcal{H}_3\mathcal{H}_4)^2=(p_2)$, so
$[\mathcal{H}_1]=[\mathcal{H}_2]$ and $[\mathcal{H}_3]=[\mathcal{H}_4]$.\\
Since
 \begin{center}
$(\mathcal{H}_0\mathcal{H}_1\mathcal{H}_3)^2=((1+i)\pi_1\pi_3)$,\\ $(\mathcal{H}_1\mathcal{H}_3)^2=(\pi_1\pi_3)$,
\end{center}
 \begin{center}
$(1+i)\pi_1\pi_3=[(eg-4fh)-2(eh+fg)]+i[(eg-4fh)+2(eh+fg)]$,\\
$\pi_1\pi_3=(eg-4fh)+2i(eh+fg)$,\\
 $[(eg-4fh)-2(eh+fg)]^2+[(eg-4fh)+2(eh+fg)]^2=2p_1p_2$ and\\
 $(eg-4fh)^2+4(eh+fg)^2=p_1p_2$,
 \end{center}
  then
  Proposition \ref{26} implies that $\mathcal{H}_0\mathcal{H}_1\mathcal{H}_3$ and $\mathcal{H}_1\mathcal{H}_3$
  are not principal in $\kk$, therefore $$\langle [\mathcal{H}_0], [\mathcal{H}_1],
  [\mathcal{H}_3]\rangle\subseteq Am_s(\kk/\QQ(i)).$$
  Proceeding as above we prove that $|Am_s(\kk/F)|=8$, which yields that
  $$Am_s(\kk/\QQ(i))=\langle [\mathcal{H}_0], [\mathcal{H}_1],
  [\mathcal{H}_3]\rangle.$$\par
 (3) Since $Q_\kk=2$,  then, from Lemma \ref{6},
 $E_\kk=\langle i, \sqrt{i\varepsilon_3}\rangle$,  $N(\varepsilon_3)=1$ and $x\pm1$ is a square in $\NN$; so, as
 $(\mathcal{H}_1\mathcal{H}_2)^2=(p_1)$  and $(\mathcal{H}_3\mathcal{H}_4)^2=(p_2)$,  Proposition \ref{28} yields that
  $[\mathcal{H}_1]\neq[\mathcal{H}_2]$ and $[\mathcal{H}_3]\neq[\mathcal{H}_4]$; moreover, as
  $(\mathcal{H}_0^2\mathcal{H}_1\mathcal{H}_2\mathcal{H}_3\mathcal{H}_4)^2=(2p_1p_2)$, then
  $\mathcal{H}_1\mathcal{H}_2\mathcal{H}_3\mathcal{H}_4=(\frac{\sqrt{2p_1p_2}}{1+i})$, hence
  $[\mathcal{H}_4]=[\mathcal{H}_1\mathcal{H}_2\mathcal{H}_3]$; therefore  Proposition \ref{26} allowed us
   to state that $$\langle [\mathcal{H}_0], [\mathcal{H}_1],
   [\mathcal{H}_2], [\mathcal{H}_3]\rangle\subseteq Am_s(\kk/\QQ(i)).$$ As above we prove that $|Am(\kk/F)|=2^4$
and $E_F\cap N_{\kk/F}(\kk^\times)=\langle i\rangle$, hence formula (\ref{50}) yields that
$|Am(\kk/F)|=|Am_s(\kk/F)|=16$, so $$Am_s(\kk/\QQ(i))=\langle [\mathcal{H}_0], [\mathcal{H}_1],
   [\mathcal{H}_2], [\mathcal{H}_3]\rangle.$$
\end{proof}
\section{\textbf{Proof of the Main Theorem}}
We know that  $\kk=\QQ(\sqrt{2p_1p_2},i)$ and its  genus field is\\ $\k=\QQ(\sqrt2, \sqrt{p_1}, \sqrt{p_2}, i)$, so $[\k:\kk]=4$ and there are three unramified quadratic  extensions of $\kk$ abelian  over $\QQ$ which are
$\KK_1=\kk(\sqrt{p_1})=\QQ(\sqrt{p_1},\sqrt{2p_2}, i)$, $\KK_2=\kk(\sqrt{p_2})=\QQ(\sqrt{p_2},\sqrt{2p_1},i)$
and $\KK_3=\kk(\sqrt{2})=\QQ(\sqrt 2, \sqrt{p_1p_2}, i)$; let  $N_i$ denote the norm  $N_{\KK_i/\kk}$. Let $\mathbf{C}_{\mathds{\kk}}$ (resp
$\mathbf{C}_{\mathds{\KK}_i}$) be the ideal class group of $\kk$ (resp $\KK_i$),
we denote by  $J_{\KK_i}$  the homomorphism  from $\mathbf{C}_{\mathds{\kk}}$ to $\mathbf{C}_{\mathds{\KK}_i}$
that maps to the class of an ideal $\mathcal{I}$ of $\kk$ the class of the ideal generated by  $\mathcal{I}$ in $\KK_i$. We keep the notations  defined in the beginning of the preceding  section. To prove the Main Theorem, we must study the capitulation problem of the 2-classes of $\kk$ in each $\KK_i$  and in $\k$.
\subsection{\textbf{Capitulation in $\KK_1$}}
Let $\varepsilon_1$ (resp $\varepsilon_2$, $\varepsilon_3$ ) denote the fundamental unit of $\QQ(\sqrt{p_1})$
(resp. $\QQ(\sqrt{2p_2})$, $\QQ(\sqrt{2p_1p_2})$). Put $\varepsilon_3=x+y\sqrt{2p_1p_2}$.
\begin{them}\label{29}
Assume that $N(\varepsilon_2)=N(\varepsilon_3)=1$, then
\begin{enumerate}[\upshape\indent(1)]
  \item If  $x\pm1$ is a square in $\NN$, then $|kerJ_{\KK_1}|=4$.
  \item Else $|kerJ_{\KK_1}|=2$.
\end{enumerate}
\end{them}
\begin{proof}
 From  Proposition $\ref{13}$,  $E_{\KK_1}=\langle i, \varepsilon_1, \sqrt{i\varepsilon_2},
 \sqrt{i\varepsilon_2\varepsilon_3}\rangle$ or $E_{\KK_1}=\langle i, \varepsilon_1, \sqrt{i\varepsilon_2},
  \sqrt{i\varepsilon_3}\rangle$, then $N_1(E_{\KK_1})=\langle i, \varepsilon_3\rangle$.\\
 \indent (1) From  Lemma \ref{6}, $E_{\kk}=\langle i,
 \sqrt{i\varepsilon_3}\rangle$, so $[E_{\kk}:N_1(E_{\KK_1})]=2$, and  the relation (\ref{1}) implies that
 $|kerJ_{\KK_1}|=4$.\\
 \indent (2) If $x\pm1$ is not a square in $\NN$, then $E_{\kk}=\langle i, \varepsilon_3\rangle$, so $[E_{\kk}:N_1(E_{\KK_1})]=1$ and $|kerJ_{\KK_1}|=2$.
\end{proof}
\begin{coro}\label{30}
We keep the  assumptions of the preceding theorem.
\begin{enumerate}[\upshape\indent(1)]
 \item If  $x\pm1$ is a square in $\NN$, then $kerJ_{\KK_1}=\langle [\mathcal{H}_1], [\mathcal{H}_2]\rangle\subset Am_s(\kk/\QQ(i))$.
  \item Else $kerJ_{\KK_1}=\langle[\mathcal{H}_1]\rangle\subset Am_s(\kk/\QQ(i))$.
\end{enumerate}
\end{coro}
\begin{proof}
As  $\mathcal{H}_1^2=(\pi_1)$, $\mathcal{H}_2^2=(\pi_2)$ and
 $\sqrt{e^2+(2f)^2}=\sqrt{p_1}\not\in\QQ(\sqrt{2p_1p_2})$, hence  Proposition \ref{26} implies that $\mathcal{H}_1$,  $\mathcal{H}_2$  are not principal in $\kk$.\\
 \indent Let us prove that $\mathcal{H}_1$, $\mathcal{H}_2$ capitulate in $\KK_1$.
 According to the equalities (\ref{8}),
 \begin{center}$\sqrt{2\pi_1\varepsilon_1}\in\KK_1$  and
 $\sqrt{2\pi_2\varepsilon_1}\in\KK_1$ or  $\sqrt{\pi_1\varepsilon_1}\in\KK_1$ and
 $\sqrt{\pi_1\varepsilon_2}\in\KK_1$,
 \end{center}
 so  putting
 \begin{center}
  $\sqrt{2\pi_1\varepsilon_1}=\alpha_1$,  $\sqrt{2\pi_2\varepsilon_1}=\alpha_2$,
  $\sqrt{\pi_1\varepsilon_1}=\beta_1$ and $\sqrt{\pi_1\varepsilon_2}=\beta_2$,
   \end{center}
   we get
 \begin{center}
    $\mathcal{H}_1^2=(\frac{\alpha_1}{1+i})^2$
  and  $\mathcal{H}_2^2=(\frac{\alpha_2}{1+i})^2$ or   $\mathcal{H}_1^2=(\beta_1)^2$ and $\mathcal{H}_2^2=(\beta_2)^2$
   \end{center}
which implies the result.\\
\indent (1) If $x\pm1$ is a square in $\NN$, then $p_1(x\pm1)$ and $2p_1(x\pm1)$
are not squares in $\NN$;  moreover  as $(p_1)=(\mathcal{H}_1\mathcal{H}_2)^2$,
then  Proposition \ref{28} yields that $\mathcal{H}_1\mathcal{H}_2$ is not principal in $\kk$, so $\mathcal{H}_1$ and $\mathcal{H}_2$ lie in  different classes. Thus
 $kerJ_{\KK_1}=\langle [\mathcal{H}_1], [\mathcal{H}_2]\rangle$ which is a subgroup of $Am_s(\kk/\QQ(i))$.\\
\indent (2) If $x\pm1$ is not a square in $\NN$, then $\mathcal{H}_1\mathcal{H}_2$ is principal  in $\kk$, so $\mathcal{H}_1$, $\mathcal{H}_2$
lie in the same class; thus $kerJ_{\KK_1}=\langle [\mathcal{H}_1]\rangle \subset Am_s(\kk/\QQ(i))$.
\end{proof}
\begin{exams}\label{48}
(1) $x\pm1$ is a square in $\NN$.\par
The first table gives  integers $d$ for which $x\pm1$ is a square in $\NN$, when the second shows that for these
integers  $\mathcal{H}_1\mathcal{H}_2$ is not principal in $\kk$ and $\mathcal{H}_1$,  $\mathcal{H}_2$ capitulate
in $\KK_1$.
\begin{longtable}{| c | c | c | c |  }
\hline
$d = 2.p_1.p_2$ & $\varepsilon_d = x+y\sqrt{d}$ & $x+1$ & $ x-1 $  \\
\hline
\endfirsthead
\hline
$d = 2.p_1.p_2$ & $\varepsilon_d = x+y\sqrt{d}$ & $x+1$ & $ x-1 $  \\
\hline
\endhead
$1394=2.41.17$ & $12545+336\sqrt{1394}$ & $12546$ & $12544 $ \\
\hline
$3298=2.97.17$ & $161603+2814\sqrt{3298}$ & $161604$ & $161602 $ \\
\hline
$15266=2.449.17$ & $1236545+10008\sqrt{15266}$ & $1236546$ & $1236544 $ \\
\hline
\end{longtable}
\begin{longtable}{| c | c | c | c | }
\hline
$d = 2.p_1.p_2$ & $\mathcal{H}_1\mathcal{H}_2$ in $\kk$ &
                             $\mathcal{H}_1$ in $\KK_1$ &  $\mathcal{H}_2$ in $\KK_1$\\
\hline
\endfirsthead
\hline
$d = 2.p_1.p_2$ &  $\mathcal{H}_1\mathcal{H}_2$ in $\kk$ &
                             $\mathcal{H}_1$ in $\KK_1$ &  $\mathcal{H}_2$ in  $\KK_1$\\
\hline
\endhead
$1394=2.41.17$ & $[0, 0, 1, 0]~$ & $[0, 0, 0, 0]~$ & $[0, 0, 0, 0]~$\\
\hline
$3298=2.97.17$ & $[4, 2, 1, 0]~$ & $[0, 0, 0, 0]~$ & $[0, 0, 0, 0]~$\\
\hline
$15266=2.449.17$ & $[0, 0, 1, 0]~$ & $[0, 0, 0, 0]~$ & $[0, 0, 0, 0]~$\\
\hline
\end{longtable}
(2) $x+1$ and $x-1$ are not squares in $\NN$.\par
In this table we give  integers $d$ for which $x+1$ and $x-1$ are not squares in $\NN$, we note that
$\mathcal{H}_1\mathcal{H}_2$ is principal in $\kk$ and   $\mathcal{H}_1$ capitulates in $\KK_1$.
\begin{longtable}{| c | c | c | c |c|c| }
\hline
$d = 2.p_1.p_2$ & $\varepsilon_d = x+y\sqrt{d}$ & $x+1$ & $ x-1 $ & $ \mathcal{H}_1\mathcal{H}_2$ in $\kk$ & $ \mathcal{H}_1$ in $\KK_1$\\
\hline
\endfirsthead
\hline
$d = 2.p_1.p_2$ & $\varepsilon_d = x+y\sqrt{d}$ & $x+1$ & $ x-1 $ & $ \mathcal{H}_1\mathcal{H}_2$ in $\kk$ & $ \mathcal{H}_1$ in $\KK_1$\\
\hline
\endhead
$890=2.5.89$ & $179+6\sqrt{890}$ & $180$ & $178 $ & $[0, 0, 0]~$&$[0, 0, 0, 0]~$\\
\hline
$1802=2.53.17$ & $849+20\sqrt{1802}$ & $850$ & $848 $ & $[0, 0, 0]~$&$[0, 0, 0, 0]~$\\
\hline
$5402=2.37.73$ & $147+2\sqrt{5402}$ & $148$ & $146 $ & $[0, 0, 0]~$&$[0, 0, 0]~$\\
\hline
\end{longtable}
\end{exams}
\begin{them}\label{31}
\begin{enumerate}[\upshape\indent(1)]
  \item If $N(\varepsilon_2)=N(\varepsilon_3)=-1$ or $N(\varepsilon_2)=-N(\varepsilon_3)=1$, then $|kerJ_{\KK_1}|=4$.
  \item Assume that $N(\varepsilon_3)=-N(\varepsilon_2)=1$, then\\
   (i) If  $x\pm1$ is a square in $\NN$, then  $|kerJ_{\KK_1}|=8$.\\
   (ii) Else $|kerJ_{\KK_1}|=4$.
\end{enumerate}
\end{them}
\begin{proof}
\indent (1) From   Propositions \ref{7} and \ref{11},    $N_1(E_{\KK_1})=\langle -1, \varepsilon_3\rangle$
or  $\langle -1, i\varepsilon_3\rangle$ or $\langle i, \varepsilon_3^2\rangle$;
on the other hand,   Lemma \ref{6} yields that  $E_{\kk}=\langle i, \varepsilon_3\rangle$, consequently
 $[E_{\kk}:N_1(E_{\KK_1})]=2$ and $|kerJ_{\KK_1}|=4$.\\
\indent (2) From Proposition \ref{12}, we get:\par
  (i) If  $x\pm1$ is a square in $\NN$, then $N_1(E_{\KK_1})=\langle -1, i\varepsilon_3\rangle$,
 as Lemma \ref{6} yields that $E_{\kk}=\langle i, \sqrt{i\varepsilon_3}\rangle$, so
  $[E_{\kk}:N_1(E_{\KK_1})]=4$ and $|kerJ_{\KK_1}|=8$.\par
 (ii) If $x\pm1$ is not a square in $\NN$, then
 $N_1(E_{\KK_1})=\langle -1, \varepsilon_3\rangle$ or
 $\langle -1, i\varepsilon_3\rangle$;
  since in this case $E_{\kk}=\langle i, \varepsilon_3\rangle$,  hence  $[E_{\kk}:N_1(E_{\KK_1})]=2$ and  $|kerJ_{\KK_1}|=4$.
\end{proof}
\begin{coro}\label{32}
\begin{enumerate}[\upshape\indent(1)]
  \item If $N(\varepsilon_2)=N(\varepsilon_3)=-1$ or $N(\varepsilon_2)=-N(\varepsilon_3)=1$,\\
  then $kerJ_{\KK_1}=\langle [\mathcal{H}_1], [\mathcal{H}_2]\rangle \subset Am_s(\kk/\QQ(i))$.
  \item Assume that $N(\varepsilon_3)=-N(\varepsilon_2)=1$, then
  \begin{enumerate}[\upshape\indent(i)]
    \item  If  $x\pm1$ is a square in $\NN$, then
   $kerJ_{\KK_1}=\langle [\mathcal{H}_1], [\mathcal{H}_2], [\mathcal{H}_0\mathcal{H}_3]\rangle \subset Am_s(\kk/\QQ(i))$.
    \item Else $ kerJ_{\KK_1}=\langle [\mathcal{H}_1], [\mathcal{H}_0\mathcal{H}_3]\rangle \subset Am_s(\kk/\QQ(i))$.
  \end{enumerate}
\end{enumerate}
\end{coro}
\begin{proof}
(1) As $N(\varepsilon_3)=-1$, then  Proposition \ref{28}  yields that $\mathcal{H}_1\mathcal{H}_2$ is not
principal in $\kk$;  even  Proposition \ref{26} implies that  $\mathcal{H}_1$,  $\mathcal{H}_2$ are not principal in $\kk$. Proceeding as in Corollary \ref{30}, we prove that $kerJ_{\KK_1}=\langle [\mathcal{H}_1],
 [\mathcal{H}_2]\rangle$. With our assumptions,  $Am_s(\kk/\QQ(i))=\langle [\mathcal{H}_0], [\mathcal{H}_1], [\mathcal{H}_2]\rangle$ (see Proposition \ref{52}), so $kerJ_{\KK_1}$ is a subgroup of $Am_s(\kk/\QQ(i))$.\par
 (2) As $(\mathcal{H}_0\mathcal{H}_1\mathcal{H}_3)^2=((1+i)\pi_1\pi_3)$, so proceeding as in Proposition \ref{52}, we get
 $\mathcal{H}_0\mathcal{H}_1\mathcal{H}_3$,
  $\mathcal{H}_0\mathcal{H}_2\mathcal{H}_3$ are not principal in $\kk$.  On the other hand, as
   $N(\varepsilon_2)=-1$, so the equalities  (\ref{9}) yields that
   $\sqrt{(1\pm i)\pi_3\varepsilon_2}\in\KK_1$, hence there exists
    $\alpha\in\KK_1$  such that $(1\pm i)\pi_3\varepsilon_2=\alpha^2$  i.e. $(\mathcal{H}_0\mathcal{H}_3)^2=(\alpha^2)$, therefore
     $\mathcal{H}_0\mathcal{H}_3$ capitulates  in $\KK_1$.\\
\indent (i) If $x\pm1$ is a square in $\NN$, then Propositions \ref{26} and \ref{28} state that $\mathcal{H}_1$,
$\mathcal{H}_2$ and $\mathcal{H}_1\mathcal{H}_2$ are not principal in $\kk$. Hence $kerJ_{\KK_1}=\langle [\mathcal{H}_1], [\mathcal{H}_2], [\mathcal{H}_0\mathcal{H}_3]\rangle $. By Proposition \ref{52} (2), we get $kerJ_{\KK_1}\subset Am_s(\kk/\QQ(i))$.\par
 (ii) If  $x\pm1$ is not a square in $\NN$, then $\mathcal{H}_1$,  $\mathcal{H}_2$ lie in the same class; so $ kerJ_{\KK_1}=\langle [\mathcal{H}_1], [\mathcal{H}_0\mathcal{H}_3]\rangle$, and by Proposition \ref{52} (1).(ii), we get $kerJ_{\KK_1}\subset Am_s(\kk/\QQ(i))$.
\end{proof}
\begin{exams}\label{49}
(1) $N(\varepsilon_2)=N(\varepsilon_3)=-1$ or $N(\varepsilon_2)=-N(\varepsilon_3)=1$.\\
This table  gives  integers $d$ for which
$\mathcal{H}_1\mathcal{H}_2$ is not principal in $\kk$ and   $\mathcal{H}_1$, $\mathcal{H}_2$  capitulate in $\KK_1$.
\begin{longtable}{| c | c | c | c |c|c| }
\hline
$d = 2.p_1.p_2$ & $N(\varepsilon_2)$ & $N(\varepsilon_3)$ & $\mathcal{H}_1\mathcal{H}_2$ in $\kk$ &
                             $\mathcal{H}_1$ in $\KK_1$ &  $\mathcal{H}_2$ in  $\KK_1$\\
\hline
\endfirsthead
\hline
$d = 2.p_1.p_2$ & $N(\varepsilon_2) $ & $ N(\varepsilon_3)$ & $\mathcal{H}_1\mathcal{H}_2$ in $\kk$ &
                             $\mathcal{H}_1$ in $\KK_1$ &  $\mathcal{H}_2$ in  $\KK_1$\\
\hline
\endhead
$290=2.5.29$ & $-1$ & $-1$ & $[5, 0, 0]~$ & $[0, 0, 0]~$ & $[0, 0, 0]~$\\
\hline
$442=2.13.17$ & $1$ & $-1$ & $[2, 0, 0]~$ & $[0, 0, 0]~$ & $[0, 0, 0]~$\\
\hline
$754=2.29.13$ & $1$ & $-1$ & $[0, 0, 1]~$ & $[0, 0, 0]~$ & $[0, 0, 0]~$\\
\hline
$1066=2.41.13$ & $-1$ & $-1$ & $[0, 1, 0]~$ & $[0, 0, 0]~$ & $[0, 0, 0]~$\\
\hline
\end{longtable}
 (2) $N(\varepsilon_3)=-N(\varepsilon_2)=1$.\par
(i) $x\pm1$ is a square in $\NN$.\par
 The first table gives  integers $d$ for which $x\pm1$ is a square in $\NN$ and $\mathcal{H}_1\mathcal{H}_2$ is
 not principal in $\kk$, when the second shows that for these
integers   $\mathcal{H}_1$,  $\mathcal{H}_2$ and  $\mathcal{H}_0\mathcal{H}_3$ capitulate
in $\KK_1$, but  $\mathcal{H}_0$,  $\mathcal{H}_3$ do not.
\small
\begin{longtable}{| c | c | c | c |c|c| }
\hline
$d = 2.p_1.p_2$ & $\varepsilon_d = x+y\sqrt{d}$ & $x+1$ & $ x-1 $ & $\mathcal{H}_1\mathcal{H}_2$ in $\kk$\\
\hline
\endfirsthead
\hline
$d = 2.p_1.p_2$ & $ \varepsilon_d = x+y\sqrt{d}$ & $x+1$ & $ x-1 $ & $\mathcal{H}_1\mathcal{H}_2$ in $\kk$\\
\hline
\endhead
$1394=2.17.41$ & $12545+336\sqrt{1394}$ & $12546$ & $12544$ & $ [0, 0, 1, 0]~$\\
\hline
$7298=2.89.41$ & $357603+4186\sqrt{7298}$ & $357604$ & $357602$ & $ [0, 0, 0, 1]~$\\
\hline
$16498=2.73.113$ & $1336337+10404\sqrt{16498}$ & $1336338$ & $1336336$ & $ [48, 0, 1, 1]~$\\
\hline
\end{longtable}
\begin{longtable}{| c | c | c | c | c | c | }
\hline
$d = 2.p_1.p_2$ & $\mathcal{H}_1$ in $\KK_1$ &
                  $\mathcal{H}_2$ in $\KK_1$ &  $\mathcal{H}_0\mathcal{H}_3$ in  $\KK_1$ &
				  $\mathcal{H}_3$ in  $\KK_1$ &  $\mathcal{H}_0$ in  $\KK_1$ \\
\hline
\endfirsthead
\hline
$d = 2.p_1.p_2$ & $\mathcal{H}_1$ in $\KK_1$ &
                  $\mathcal{H}_2$ in $\KK_1$ &  $\mathcal{H}_0\mathcal{H}_3$ in  $\KK_1$ &
				  $\mathcal{H}_3$ in  $\KK_1$ &  $\mathcal{H}_0$ in  $\KK_1$\\
\hline
\endhead
$1394=2.17.41$ & $[0, 0, 0, 0]~$ & $[0, 0, 0, 0]~$ & $[0, 0, 0, 0]~$ & $[24, 0, 0, 0]~$ & $[24, 0, 0, 0]~$\\
\hline
$7298=2.89.41$ & $[0, 0, 0, 0]~$ & $[0, 0, 0, 0]~$ & $[0, 0, 0, 0]~$ & $[42, 2, 2, 0]~$ & $[42, 2, 2, 0]~$\\
\hline
$16498=2.73.113$ & $[0, 0, 0, 0]~$ & $[0, 0, 0, 0]~$ & $[0, 0, 0, 0]~$ & $[96, 0, 0, 0]~$ & $[96, 0, 0, 0]~$\\
\hline
\end{longtable}
\normalsize
(ii) $x+1$ and $x-1$ are not  squares in $\NN$.\par
 The first table gives  integers $d$ for which $x+1$ and $x-1$ are not  squares in $\NN$ and $\mathcal{H}_1\mathcal{H}_2$ is
 principal in $\kk$, when the second shows that for these
integers   $\mathcal{H}_1$ and  $\mathcal{H}_0\mathcal{H}_3$ capitulate
in $\KK_1$, but  $\mathcal{H}_0$,  $\mathcal{H}_3$ do not.
\footnotesize
\begin{longtable}{| c | c | c | c |c|c| }
\hline
$d = 2.p_1.p_2$ & $\varepsilon_d = x+y\sqrt{d}$ & $x+1$ & $ x-1 $ & $\mathcal{H}_1\mathcal{H}_2$ in $\kk$\\
\hline
\endfirsthead
\hline
$d = 2.p_1.p_2$ & $ \varepsilon_d = x+y\sqrt{d}$ & $x+1$ & $ x-1 $ & $\mathcal{H}_1\mathcal{H}_2$ in $\kk$\\
\hline
\endhead
$410=2.5.41 $ & $81+4\sqrt{410}$ & $82$ & $80$ & $ [0, 0, 0]~$\\
\hline
$2938=2.13.113 $ & $786707+14514\sqrt{2938}$ & $786708$ & $786706$ & $ [0, 0, 0]~$\\
\hline
$3034=2.37.41 $ & $4055973299+73635510\sqrt{3034}$ & $4055973300$ & $4055973298$ & $ [0, 0, 0]~$\\
\hline
$8090=2.5.809 $ & $1619+18\sqrt{8090}$ & $1620$ & $1618$ & $ [0, 0, 0]~$\\
\hline
\end{longtable}
\normalsize
\begin{longtable}{| c | c | c | c | c |  }
\hline
$d = 2.p_1.p_2$ & $\mathcal{H}_1$ in $\KK_1$  &  $\mathcal{H}_0\mathcal{H}_3$ in  $\KK_1$ &
				  $\mathcal{H}_3$ in  $\KK_1$ &  $\mathcal{H}_0$ in  $\KK_1$ \\ \hline
\endfirsthead \hline
$d = 2.p_1.p_2$ & $\mathcal{H}_1$ in $\KK_1$ &   $\mathcal{H}_0\mathcal{H}_3$ in  $\KK_1$ &
				  $\mathcal{H}_3$ in  $\KK_1$ &  $\mathcal{H}_0$ in  $\KK_1$\\ \hline \endhead
$410=2.5.41 $ & $[0, 0, 0]~$ & $[0, 0, 0]~$ & $[8, 0, 0]~$ & $[8, 0, 0]~$\\
\hline
$2938=2.13.113 $ & $[0, 0, 0]~$ & $[0, 0, 0]~$  & $[24, 0, 0]~$ & $[24, 0, 0]~$\\
\hline
$3034=2.37.41 $ & $[0, 0, 0]~$ & $[0, 0, 0]~$  & $[0, 0, 2]~$ & $[0, 0, 2]~$\\
\hline
$8090=2.5.809 $ & $[0, 0, 0, 0]~$ & $[0, 0, 0, 0]~$ & $[168, 0, 0, 0]~$ & $[168, 0, 0, 0]~$\\
\hline
\end{longtable}
\end{exams}
\subsection{\textbf{Capitulation in $\KK_2$}}
As $p_1$,  $p_2$  play symmetric roles, then  putting  $\varepsilon_1$ (resp $\varepsilon_2$,
$\varepsilon_3$ )
the fundamental unit of $\QQ(\sqrt{p_2})$ (resp $\QQ(\sqrt{2p_1})$, $\QQ(\sqrt{2p_1p_2})$) and proceeding as above, we get the following results. Put $\varepsilon_3=x+y\sqrt{2p_1p_2}$.
\begin{them}\label{33}
Assume that $N(\varepsilon_2)=N(\varepsilon_3)= 1$.
\begin{enumerate}[\upshape\indent(1)]
  \item If  $x\pm1$ is a square in $\NN$, then $|kerJ_{\KK_2}|=4$.
  \item Else $|kerJ_{\KK_2}|=2$.
\end{enumerate}
\end{them}
\begin{coro}\label{34}
We keep the  assumptions of the Theorem $\ref{33}$, then
\begin{enumerate}[\upshape\indent(1)]
 \item If  $x\pm1$ is a square in $\NN$, then $kerJ_{\KK_2}=\langle [\mathcal{H}_3], [\mathcal{H}_4]\rangle \subset Am_s(\kk/\QQ(i))$.
  \item Else $kerJ_{\KK_2}=\langle[\mathcal{H}_3]\rangle \subset Am_s(\kk/\QQ(i))$.
\end{enumerate}
\end{coro}
\begin{them}\label{35}
\begin{enumerate}[\upshape\indent(1)]
  \item If $N(\varepsilon_2)=N(\varepsilon_3)=-1$ or $N(\varepsilon_2)=-N(\varepsilon_3)=1$, then $|kerJ_{\KK_2}|=4$.
  \item If $N(\varepsilon_3)=-N(\varepsilon_2)=1$, then\\
   (i) If $x\pm1$ is a square in $\NN$, then  $|kerJ_{\KK_2}|=8$.\\
   (ii) Else $|kerJ_{\KK_2}|=4$.
\end{enumerate}
\end{them}
\begin{coro}\label{36}
\begin{enumerate}[\upshape\indent(1)]
  \item If $N(\varepsilon_2)=N(\varepsilon_3)=-1$ or $N(\varepsilon_2)=-N(\varepsilon_3)=1$,\\
  then $kerJ_{\KK_2}=\langle [\mathcal{H}_3], [\mathcal{H}_4]\rangle \subset Am_s(\kk/\QQ(i))$.
  \item Assume that $N(\varepsilon_3)=-N(\varepsilon_2)=1$,  then  \\
   (i) If  $x\pm1$ is a square in $\NN$, then\\
   $kerJ_{\KK_2}=\langle[\mathcal{H}_0\mathcal{H}_1],  [\mathcal{H}_3], [\mathcal{H}_4]\rangle \subset Am_s(\kk/\QQ(i))$.\\
   (ii) Else $kerJ_{\KK_2}=\langle[\mathcal{H}_0\mathcal{H}_1], [\mathcal{H}_3]\rangle \subset Am_s(\kk/\QQ(i))$.
\end{enumerate}
\end{coro}
\subsection{\textbf{Capitulation in $\KK_3$}}
Let $\varepsilon_1$ (resp. $\varepsilon_2$, $\varepsilon_3$ ) denote the  fundamental unit of
$\QQ(\sqrt{2})$ (resp. $\QQ(\sqrt{p_1p_2})$, $\QQ(\sqrt{2p_1p_2})$) and $q=q(\KK_3^+/\QQ)$.
 \begin{them}\label{37}
Suppose   $N(\varepsilon_2)=N(\varepsilon_3)=-1$ or  $N(\varepsilon_2)=-N(\varepsilon_3)=1$, then
\begin{enumerate}[\upshape\indent(1)]
  \item If  $q=1$, then $|kerJ_{\KK_3}|=4$.
  \item If  $q=2$, then $|kerJ_{\KK_3}|=2$.
\end{enumerate}
\end{them}
\begin{proof} Note first that $\sqrt i=\frac{1+i}{\sqrt2}\in\KK_3$.\\
\indent (1) From  Propositions \ref{15} and \ref{16},  we get
$E_{\KK_3}=\langle \sqrt i,\varepsilon_1, \varepsilon_2,  \varepsilon_3\rangle$, so
$N_3(E_{\KK_3})=\langle i, \varepsilon_3^2\rangle$;  as, from Lemma \ref{6},  $E_{\kk}=\langle i, \varepsilon_3\rangle$, then
 $[E_{\kk}:N_3(E_{\KK_3})]=2$ and  $|kerJ_{\KK_3}|=4$.\par
 (2)  Proposition \ref{16} yields that
$E_{\KK_3}=\langle \sqrt i,\varepsilon_1, \varepsilon_2,  \sqrt{\varepsilon_1 \varepsilon_2\varepsilon_3}\rangle$, so
$N_3(E_{\KK_3})=\langle i, \varepsilon_3\rangle$; on the other hand,   Lemma \ref{6} implies that
$E_{\kk}=\langle i, \varepsilon_3\rangle$, hence $[E_{\kk}:N_3(E_{\KK_3})]=1$ and
  $|kerJ_{\KK_3}|=2$.
\end{proof}
\begin{coro}\label{38}
Keep the assumptions of Theorem $\ref{37}$, then
\begin{enumerate}[\upshape\indent(1)]
  \item If  $q=2$, then $kerJ_{\KK_3}=\langle[\mathcal{H}_0]\rangle \subset Am_s(\kk/\QQ(i))$.
  \item If  $q=1$, then $kerJ_{\KK_3}=\langle[\mathcal{H}_0], [\mathcal{H}_1\mathcal{H}_2]\rangle \subset Am_s(\kk/\QQ(i))$.
\end{enumerate}
\end{coro}
\begin{proof}
Note that
$\sqrt{(1+i)\varepsilon_1}=\frac{1}{2}(2+(1+i)\sqrt2)$, so
 there  exists $\beta\in\KK_3$ such that $\mathcal{H}_0^2=(1+i)=(\beta^2)$, hence $\mathcal{H}_0$
capitulates in $\KK_3$. \\
 \indent (1) If $q=2$, then $kerJ_{\KK_3}=\langle[\mathcal{H}_0]\rangle$, hence by Proposition \ref{52} (1).(i),  $kerJ_{\KK_3} \subset Am_s(\kk/\QQ(i))$.\\
\indent (2) There are two cases to distinguish.\\
 \indent (i) If  $N(\varepsilon_2)=1$, then, from the proof of Proposition
\ref{15}, we get that $\sqrt{p_1\varepsilon_2}\in\KK_3$, so  putting  $\alpha=\sqrt{p_1\varepsilon_3}$,
 we deduce  that $(\alpha^2)=\mathcal{H}_1^2\mathcal{H}_2^2$, which implies that $\mathcal{H}_1\mathcal{H}_2$
 capitulates in $\KK_3$.\\
 \indent (ii) If   $N(\varepsilon_2)=-1$, then $\varepsilon_1\varepsilon_2\varepsilon_3$
 is not a square in $\QQ(\sqrt2, \sqrt{p_1p_2})$,  so, according to the remark  \ref{17}, we get
 $\sqrt{p_1\varepsilon_1\varepsilon_2\varepsilon_3}\in\KK_3$, hence $\mathcal{H}_1\mathcal{H}_2$
 capitulates in $\KK_3$.\\
  \indent On the other hand,   Proposition \ref{28}  yields that $\mathcal{H}_1\mathcal{H}_2$
 is not principal in $\kk$, furthermore  $(\mathcal{H}_0\mathcal{H}_1\mathcal{H}_2)^2=((1+i)p_1)$ and
 $\sqrt2\not\in\kk$, then  Proposition \ref{26} implies  that $\mathcal{H}_0\mathcal{H}_1\mathcal{H}_2$ is not principal in
  $\kk$ i.e. $\mathcal{H}_0$ and $\mathcal{H}_1\mathcal{H}_2$ lie in different classes; hence the result. By  Proposition \ref{52} (1).(i) we get  $kerJ_{\KK_3} \subset Am_s(\kk/\QQ(i))$.
\end{proof}
\begin{exams}
(1) First case $q=2$.
\begin{longtable}{| c | c | c | c | c | }
\hline
$d = 2.p_1.p_2$ & $q$ & $N(\varepsilon_2)$ & $N(\varepsilon_3)$ &
                        $\mathcal{H}_0$ in $\KK_3$ \\
\hline
\endfirsthead
\hline
$d = 2.p_1.p_2$ & $q$ & $N(\varepsilon_2)$ & $N(\varepsilon_3)$ &
                        $\mathcal{H}_0$ in $\KK_3$\\
\hline
\endhead
$130=2.5.13$ & $2$ & $-1$ & $-1 $ & $ [0, 0]~$\\
\hline
$1066=2.13.41$ & $2$ & $-1$ & $-1 $ & $[0, 0, 0, 0]~$\\
\hline
$2146=2.29.37$ & $2$ & $-1$ & $-1 $ & $ [0, 0]~$\\
\hline
\end{longtable}
(2) Second  case $q=1$.
\begin{longtable}{| c | c | c | c | c | c | c |}
\hline
$d = 2.p_1.p_2$ & $q$ & $N(\varepsilon_2)$ & $N(\varepsilon_3)$ &
                        $\mathcal{H}_0$ in $\KK_3$ & $\mathcal{H}_1\mathcal{H}_2$ in $\KK_3$ &
                        $\mathcal{H}_1\mathcal{H}_2$ in $\kk$\\
\hline
\endfirsthead
\hline
$d = 2.p_1.p_2$ & $q$ & $N(\varepsilon_2)$ & $N(\varepsilon_3)$ &
                        $\mathcal{H}_0$ in $\KK_3$ & $\mathcal{H}_1\mathcal{H}_2$ in $\KK_3$ &
                        $\mathcal{H}_1\mathcal{H}_2$ in $\kk$\\
\hline
\endhead
$290=2.5.29$ & $1$ & $-1$ & $-1 $ & $ [0, 0]~ $ & $ [0, 0]~ $ & $ [5, 0, 0]~$\\
\hline
$754=2.13.29$ & $1$ & $1$ & $-1 $ & $ [0, 0]~ $ & $ [0, 0]~ $ & $ [0, 0, 1]~$\\
\hline
$962=2.37.13$ & $1$ & $-1$ & $-1 $ & $ [0, 0]~ $ & $ [0, 0]~ $ & $ [7, 1, 1]~$\\
\hline
$1378=2.53.13$ & $1$ & $1$ & $-1 $ & $ [0, 0]~ $ & $ [0, 0]~ $ & $ [5, 0, 0]~$\\
\hline
\end{longtable}
\end{exams}
\begin{them}\label{39}
 If  $N(\varepsilon_3)=-N(\varepsilon_2)=1$, then  $|kerJ_{\KK_3}|=4$.
\end{them}
\begin{proof}
 If $\varepsilon_3=x+y\sqrt{2p_1p_2}$,  then  from Proposition \ref{21} we get:\par
$\bullet$ If $x\pm1$ is a square in $\NN$, then
$N_3(E_{\KK_3})=\langle i, \varepsilon_3\rangle$; as $E_{\kk}=\langle i, \sqrt{i\varepsilon_3}\rangle$, so
$[E_{\kk}:N_3(E_{\KK_3})]=2$ and $|kerJ_{\KK_3}|=4$.\par
 $\bullet$ If $x\pm1$ is not a square in $\NN$, then  $N_3(E_{\KK_3})=\langle i, \varepsilon_3^2\rangle$; as $E_{\kk}=\langle i, \varepsilon_3\rangle$, hence
 $[E_{\kk}:N_3(E_{\KK_3})]=2$ and $|kerJ_{\KK_3}|=4$.
\end{proof}
\begin{coro}\label{46}
Keep the assumptions of Theorem $\ref{39}$, then
\begin{center}
 $kerJ_{\KK_3}=\langle[\mathcal{H}_0], [\mathcal{H}_1\mathcal{H}_3]\rangle  \subset Am_s(\kk/\QQ(i))$.\end{center}
\end{coro}
\begin{proof} We have proved that $\mathcal{H}_0$ capitulates in  $\KK_3$.\par
  $\mathcal{H}_1\mathcal{H}_3$ capitulates in $\KK_3$
if and only if there exist a unit $\varepsilon\in\KK_3$, $\alpha\in\KK_3$ such that   $$\alpha^2=\varepsilon\pi_1\pi_3,$$
 as $\pi_1\pi_3=(e+2if)(g+2ih)=(eg-4fh)+2i(fg+eh)$, so putting
 $\alpha=\alpha_1+i\alpha_2$, where $\alpha_j\in\QQ(\sqrt2, \sqrt{p_1p_2})$), and choosing  $\varepsilon$ real,
  we get
    $$\alpha_1^2-\alpha_2^2=\varepsilon(eg-4fh) \text{ and }
 \alpha_1\alpha_2=\varepsilon(fg+eh),$$
 hence
 $$\alpha_1^4-\varepsilon(eg-4fh)\alpha_1^2-\varepsilon^2(fg+eh)^2=0,$$
 since
 $p_1p_2=(eg-4fh)^2+4(fg+eh)^2$,
 so $$\alpha_1^2=\frac{\varepsilon}{2}[(eg-4fh)+\sqrt{p_1p_2}]=
 \frac{\varepsilon}{4}[\pi_1\pi_3+\pi_2\pi_4+2\sqrt{p_1p_2}],$$
 therefore
 $$\alpha_1=\frac{\sqrt{\varepsilon}}{2}[\sqrt{\pi_1\pi_3}+\sqrt{\pi_2\pi_4}].$$
  Then if
 $\varepsilon=\varepsilon_2$ and $\sqrt{\varepsilon_2}=y_1\sqrt{\pi_1\pi_3}+y_2\sqrt{\pi_2\pi_4}$,
  we get  $$\alpha_1=\frac{1}{2}(a+b\sqrt{p_1p_2})\text{ and }\alpha_2=\frac{\varepsilon(fg+eh)}{\alpha_1},$$
   where $a$ and $b$ are in $\QQ$; which implies that  $\mathcal{H}_1\mathcal{H}_3$ capitulates in $\KK_3$.\\
   \indent On the other hand,  $\mathcal{H}_0$ and $\mathcal{H}_1\mathcal{H}_3$
   lie in different classes, in fact
   \begin{center}  $(\mathcal{H}_0\mathcal{H}_1\mathcal{H}_3)^2=((1+i)\pi_1\pi_3)=$\\$
   [(eg-4fh)-2(eh+fg)]+i[(eg-4fh)+2(eh+fg)]$\\
   and  $[(eg-4fh)-2(eh+fg)]^2+[(eg-4fh)+2(eh+fg)]^2=2p_1p_2$,
   \end{center}
    so  Proposition \ref{26}
   yields the result. By  Proposition \ref{52} (1).(ii) and (2) we get  $kerJ_{\KK_3} \subset Am_s(\kk/\QQ(i))$.
\end{proof}
\begin{exams}
For the case: $N(\varepsilon_3)=-N(\varepsilon_2)=1$,   $\mathcal{H}_0$,
$ \mathcal{H}_1\mathcal{H}_3$ capitulate in $\KK_3$ and $ \mathcal{H}_1$, $\mathcal{H}_3$ do not.
\begin{longtable}{| c | c | c | c | c | }
\hline
$d = 2.p_1.p_2$ & $\mathcal{H}_0$ in $\KK_3$ & $\mathcal{H}_1$ in $\KK_3$ &
                           $\mathcal{H}_3$ in $\KK_3$ &  $ \mathcal{H}_1\mathcal{H}_3$ in $\KK_3$\\
\hline
\endfirsthead
\hline
$d = 2.p_1.p_2$ & $\mathcal{H}_0$ in $\KK_3$ & $\mathcal{H}_1$ in $\KK_3$ &
                           $\mathcal{H}_3$ in $\KK_3$ &  $ \mathcal{H}_1\mathcal{H}_3$ in $\KK_3$\\
\hline
\endhead
$890=2.5.89$ & $[0, 0, 0]~$ & $[6, 0, 2]~ $ & $[6, 0, 2]~ $ & $[0, 0, 0]~$\\
\hline
$1802=2.53.17$ & $[0, 0, 0]~$ & $[48, 0, 0]~ $ & $[48, 0, 0]~ $ & $[0, 0, 0]~$\\
\hline
$5002=2.61.41$ & $[0, 0, 0, 0]~$ & $[112, 0, 0, 0]~ $ & $[112, 0, 0, 0]~ $ & $[0, 0, 0, 0]~$\\
\hline
\end{longtable}
\end{exams}
\begin{them}\label{45}
 Suppose  $N(\varepsilon_2)=N(\varepsilon_3)=1$ and put $\varepsilon_3=x+y\sqrt{2p_1p_2}$, then
\begin{enumerate}[\upshape\indent(1)]
  \item If  $x\pm1$ is a square in $\NN$, then $|kerJ_{\KK_3}|=4$.
  \item Else  $|kerJ_{\KK_3}|=2$.
\end{enumerate}
\end{them}
\begin{proof}
 From  Proposition \ref{22}, we get\par
$(1)$
$N_3(E_{\KK_3})=\langle i, \varepsilon_3\rangle$; as $E_{\kk}=\langle i, \sqrt{i\varepsilon_3}\rangle$, so
$[E_{\kk}:N_3(E_{\KK_3})]=2$ and $|kerJ_{\KK_3}|=4$.\par
 $(2)$   $N_3(E_{\KK_3})=\langle i, \varepsilon_3\rangle$; since $E_{\kk}=\langle i, \varepsilon_3\rangle$, hence
 $[E_{\kk}:N_3(E_{\KK_3})]=1$ and $|kerJ_{\KK_3}|=2$.
\end{proof}
\begin{coro}\label{40}
We keep  the assumptions of Theorem $\ref{45}$, then
\begin{enumerate}[\rm\indent(1)]
\item If  $x\pm1$ is a square in $\NN$, then $kerJ_{\KK_3}=\langle[\mathcal{H}_0],
[\mathcal{H}_1\mathcal{H}_2]\rangle  \subset Am_s(\kk/\QQ(i))$.
 \item  Else $kerJ_{\KK_3}=\langle[\mathcal{H}_0]\rangle  \subset Am_s(\kk/\QQ(i))$.
\end{enumerate}
\end{coro}
\begin{proof} Note first  that $\mathcal{H}_0$ capitulates in  $\KK_3$.\par
 (1)  Proposition \ref{28} states that
 $\mathcal{H}_1\mathcal{H}_2$ is not principal in $\kk$; furthermore, from Equality (\ref{23}), we get  $p_1\varepsilon_2$
 is a square in $\KK_3$, so $\mathcal{H}_1\mathcal{H}_2$ capitulates in $\KK_3$. As above we prove that  $\mathcal{H}_0$, $\mathcal{H}_1\mathcal{H}_3$ lie in different classes, so $kerJ_{\KK_3}=\langle[\mathcal{H}_0], [\mathcal{H}_1\mathcal{H}_2]\rangle$, and from proposition \ref{52},  $kerJ_{\KK_3} \subset Am_s(\kk/\QQ(i))$.\par
  (2) As $\mathcal{H}_0$ capitulates in $\KK_3$, so
 $kerJ_{\KK_3}=\langle[\mathcal{H}_0]\rangle$ which is a subgroup of $ Am_s(\kk/\QQ(i))$.
\end{proof}
\begin{exams}
(1) $x\pm1$ is a square in $\NN$.
 The first table gives  integers $d$ for which $x\pm1$ is square in $\NN$ and $\mathcal{H}_1\mathcal{H}_2$ is
 not principal in $\kk$, when the second shows that for these
integers   $\mathcal{H}_0$ and  $\mathcal{H}_1\mathcal{H}_2$ capitulate
in $\KK_3$, but  $\mathcal{H}_1$,  $\mathcal{H}_2$ do not.
\small
\begin{longtable}{| c | c | c | c | c | }
\hline
$d = 2.p_1.p_2$ & $\varepsilon_d = x+y\sqrt{d}$ & $x+1$ & $ x-1 $ & $\mathcal{H}_1\mathcal{H}_2$ in $\kk$\\
\hline
\endfirsthead
\hline
$d = 2.p_1.p_2$ & $ \varepsilon_d = x+y\sqrt{d}$ & $x+1$ & $ x-1 $ & $\mathcal{H}_1\mathcal{H}_2$ in $\kk$\\
\hline
\endhead
$12994=2.73.89$ & $12995+114\sqrt{12994 }$ & $12996 $ & $12994 $ & $[0, 6, 2, 0]~$\\
\hline
$14722=2.17.433$ & $132497+1092\sqrt{14722 }$ & $132498 $ & $132496 $ & $[8, 0, 0, 0]~$\\
\hline
$32882=2.41.401$ & $295937+1632\sqrt{32882 }$ & $295938 $ & $295936 $ & $[48, 0, 0, 0]~$\\
\hline
$46658=2.41.569$ & $46657+216\sqrt{46658 }$ & $46658 $ & $46656 $ & $[64, 0, 0, 0]~$\\
\hline
\end{longtable}
\tiny
\begin{longtable}{| c | c | c | c | c | c |}
\hline
$d = 2.p_1.p_2$ & $\mathcal{H}_0$ in $\KK_3$ & $\mathcal{H}_1$ in $\KK_3$ &
                           $\mathcal{H}_2$ in $\KK_3$ &  $ \mathcal{H}_1\mathcal{H}_2$ in $\KK_3$\\
\hline
\endfirsthead
\hline
$d = 2.p_1.p_2$ & $\mathcal{H}_0$ in $\KK_3$ & $\mathcal{H}_1$ in $\KK_3$ &
                           $\mathcal{H}_1$ in $\KK_3$ &  $ \mathcal{H}_1\mathcal{H}_2$ in $\KK_3$\\
\hline
\endhead
$12994=2.73.89$ & $[0, 0, 0, 0, 0, 0]~$ & $[0, 0, 0, 1, 1, 1]~ $ & $[0, 0, 0, 1, 1, 1]~ $ & $[0, 0, 0, 0, 0, 0]~ $ \\
\hline
$14722=2.17.433$ & $[0, 0, 0, 0, 0, 0]~$ & $[56, 4, 0, 0, 0, 0]~ $ & $[56, 4, 0, 0, 0, 0]~ $ & $[0, 0, 0, 0, 0, 0]~ $ \\
\hline
$32882=2.41.401$ & $[0, 0, 0, 0, 0]~$ & $[0, 4, 0, 0, 0]~ $ & $[0, 4, 0, 0, 0]~ $ & $[0, 0, 0, 0, 0]~ $ \\
\hline
$46658=2.41.569$ & $[0, 0, 0, 0, 0, 0]~$ & $[480, 8, 0, 0, 0, 0]~ $ & $[480, 8, 0, 0, 0, 0]~ $ & $[0, 0, 0, 0, 0, 0]~ $ \\
\hline
\end{longtable}
\normalsize
(2) $x+1$ and $x-1$ are not squares in $\NN$, $\mathcal{H}_0$ capitulates in $\KK_3$.
\begin{longtable}{| c | c | c | c | c | }
\hline
$d = 2.p_1.p_2$ & $\varepsilon_d = x+y\sqrt{d}$ & $x+1$ & $ x-1 $ & $\mathcal{H}_0$ in $\KK_3$\\
\hline
\endfirsthead
\hline
$d = 2.p_1.p_2$ & $ \varepsilon_d = x+y\sqrt{d}$ & $x+1$ & $ x-1 $ & $\mathcal{H}_0$ in $\KK_3$\\
\hline
\endhead
$410=2.5.41$ & $81+4\sqrt{410 }$ & $82 $ & $80 $ & $[0, 0, 0, 0]~$\\
\hline
$2938=2.13.113$ & $786707+14514\sqrt{2938 }$ & $786708 $ & $786706 $ & $[0, 0, 0, 0]~$\\
\hline
$4010=2.401.5$ & $7219+114\sqrt{4010 }$ & $7220 $ & $7218 $ & $[0, 0, 0]~$\\
\hline
$5402=2.37.73$ & $147+2\sqrt{5402 }$ & $148 $ & $146 $ & $[0, 0, 0]~$\\
\hline
\end{longtable}
\end{exams}
\subsection{\textbf{Capitulation in $\k$}}
\begin{them}\label{47}
Let  $p_1\equiv p_2\equiv1\pmod4$ be primes and
 $\k$ the genus field of $\kk$. Let
 $Q$ denote the unit index  of the field $\QQ(\sqrt2, \sqrt{p_1p_2})$ and $\varepsilon_3$ the fundamental unit of $\QQ(\sqrt{2p_1p_2})$.
\begin{enumerate}[\upshape\indent(1)]
  \item If $N(\varepsilon_3)= -1$, then $Am_s(\kk/\QQ(i))=\langle\mathcal{H}_0, \mathcal{H}_1, \mathcal{H}_2\rangle\subseteq kerJ_{\k}$.
  \item If $N(\varepsilon_3)=1$, then
  \begin{enumerate}[\upshape\indent(i)]
    \item If $Q=2$, then $Am_s(\kk/\QQ(i))=\langle\mathcal{H}_0, \mathcal{H}_1, \mathcal{H}_2,\mathcal{H}_3 \rangle\subseteq kerJ_{\k}$.
    \item If $Q=1$, then $Am_s(\kk/\QQ(i))=\langle\mathcal{H}_0, \mathcal{H}_1, \mathcal{H}_3\rangle\subseteq kerJ_{\k}$.
  \end{enumerate}
\end{enumerate}
\end{them}
\begin{proof}
 (1)  Propositions
\ref{26},  \ref{28} imply that the classes of $\mathcal{H}_i$, where $i\in\{0, 1, 2, 3, 4\}$, are  pairwise  different. On the other hand, from Corollaries
   \ref{32}, \ref{36},  \ref{38} and \ref{40} we infer that  $[\mathcal{H}_i]\in kerJ_{\k}$; furthermore   Lemma \ref{42} states that $[\mathcal{H}_3]=[\mathcal{H}_0 \mathcal{H}_1]$ and
     $[\mathcal{H}_4]=[\mathcal{H}_0 \mathcal{H}_2]$ or $[\mathcal{H}_4]=[\mathcal{H}_0 \mathcal{H}_1]$ and
      $[\mathcal{H}_3]=[\mathcal{H}_0 \mathcal{H}_2]$, hence Proposition \ref{52} yields that
      $Am_s(\kk/\QQ(i))=\langle\mathcal{H}_0, \mathcal{H}_1, \mathcal{H}_2\rangle\subseteq kerJ_{\k}$.\par
 (2) (i) If  $Q=2$, then  Corollaries \ref{30},  \ref{32},  \ref{34},  \ref{36}  and  \ref{40}
imply that \\$\langle\mathcal{H}_0, \mathcal{H}_1, \mathcal{H}_2,\mathcal{H}_3, \mathcal{H}_4\rangle\subseteq kerJ_{\k}$.
Moreover   Propositions \ref{26},  \ref{28} state that the ideals
$\mathcal{H}_1 \mathcal{H}_2\mathcal{H}_3$, $\mathcal{H}_i \mathcal{H}_j$ and
$\mathcal{H}_0\mathcal{H}_i \mathcal{H}_j$ are not principal in $\kk$ for all
$i\in\{1, 2\}$ and $j\in\{3, 4\}$; on the other hand,
$(\mathcal{H}_0^2 \mathcal{H}_1 \mathcal{H}_2\mathcal{H}_3 \mathcal{H}_4)^2=(d)$,
so $\mathcal{H}_1 \mathcal{H}_2\mathcal{H}_3 \mathcal{H}_4=(\frac{\sqrt d}{1+i})$,
hence $[\mathcal{H}_4]=[\mathcal{H}_1 \mathcal{H}_2\mathcal{H}_3]$, and the result derived.\par
(ii) If $Q=1$, then  $\mathcal{H}_1$ and $\mathcal{H}_2$ (resp. $\mathcal{H}_3$ and $\mathcal{H}_4$)
lie in  the same class, so from Corollaries  \ref{30},  \ref{32},  \ref{36} and  \ref{38}, we get
$Am_s(\kk/\QQ(i))=\langle\mathcal{H}_0, \mathcal{H}_1, \mathcal{H}_3\rangle\subseteq kerJ_{\k}$.\\
 And this Completes the demonstration of the Main Theorem.
\end{proof}
\begin{exams}
(1) First case N($\varepsilon_d)=-1$.
\begin{longtable}{| c | c | c | c | c | }
\hline
$d = 2.p_1.p_2$ & $N(\varepsilon_d)$ & $\mathcal{H}_0$ in $\k$ & $\mathcal{H}_1$ in $\k$ &
                             $\mathcal{H}_2$ in $\k$\\
\hline
\endfirsthead
\hline
$d = 2.p_1.p_2$ & $N(\varepsilon_d)$ & $\mathcal{H}_0$ in $\k$ & $\mathcal{H}_1$ in $\k$ &
                             $\mathcal{H}_2$ in $\k$\\
\hline
\endhead
$442=2.13.17$ & $-1$ & $[0, 0, 0, 0, 0]~$ & $[0, 0, 0, 0, 0]~$ & $[0, 0, 0, 0, 0]~$\\
\hline
$1066=2.41.13$ & $-1$ & $[0, 0, 0]~$ & $[0, 0, 0]~$ & $[0, 0, 0]~$ \\
\hline
$1258=2.17.37$ & $-1$ & $[0, 0]~$ & $[0, 0]~$ & $[0, 0]~$ \\
\hline
\end{longtable}
(2) Second case $N(\varepsilon_d)=1$.\\
(i)  $Q=2$.
The first table gives  examples where $\mathcal{H}_1\mathcal{H}_2$ and $\mathcal{H}_3\mathcal{H}_4$ are
 not principal in $\kk$, when the second shows that $\mathcal{H}_0$,  $\mathcal{H}_1$, $\mathcal{H}_2$ and
 $\mathcal{H}_3$ capitulate in $\k$.
\begin{longtable}{| c | c | c | c | c |}
\hline
$d = 2.p_1.p_2$ & $N(\varepsilon_d)$ & $Q$ & $\mathcal{H}_1\mathcal{H}_2$ in $\kk$ & $\mathcal{H}_3\mathcal{H}_4$ in $\kk$\\
\hline
\endfirsthead
\hline
$d = 2.p_1.p_2$ & $N(\varepsilon_d)$ & $Q$ & $\mathcal{H}_1\mathcal{H}_2$ in $\kk$ & $\mathcal{H}_3\mathcal{H}_4$ in $\kk$\\
\hline
\endhead
$1394=2.17.41$ & $1$ & $2$ & $[0, 0, 1, 0]~$ & $[0, 0, 1, 0]~$ \\
\hline
$3298=2.97.17$ & $1$ & $2$ & $[4, 2, 1, 0]~$ & $[4, 2, 1, 0]~$ \\
\hline
$3842=2.17.113$ & $1$ & $2$ & $[0, 0, 1, 0]~$ & $[0, 0, 1, 0]~$ \\
\hline
\end{longtable}
\scriptsize
\begin{longtable}{| c | c | c | c | c | }
\hline
$d = 2.p_1.p_2$ &  $\mathcal{H}_0$ in $\k$ & $\mathcal{H}_1$ in $\k$ & $\mathcal{H}_2$ in $\k$ & $\mathcal{H}_3$ in $\k$\\
\hline
\endfirsthead
\hline
$d = 2.p_1.p_2$ & $\mathcal{H}_0$ in $\k$ & $\mathcal{H}_1$ in $\k$ & $\mathcal{H}_2$ in $\k$ & $\mathcal{H}_3$ in $\k$\\
\hline
\endhead
$1394=2.17.41$ &  $[0, 0, 0, 0, 0, 0]~$ & $[0, 0, 0, 0, 0, 0]~$ & $[0, 0, 0, 0, 0, 0]~$ & $[0, 0, 0, 0, 0, 0]~$ \\
\hline
$3298=2.97.17$ & $[0, 0, 0, 0, 0, 0, 0]~$ & $[0, 0, 0, 0, 0, 0, 0]~$ & $[0, 0, 0, 0, 0, 0, 0]~$ & $[0, 0, 0, 0, 0, 0, 0]~$\\
\hline
$3842=2.17.113$  & $[0, 0, 0, 0, 0, 0, 0]~$ & $[0, 0, 0, 0, 0, 0, 0]~$ & $[0, 0, 0, 0, 0, 0, 0]~$ & $[0, 0, 0, 0, 0, 0, 0]~$ \\
\hline
\end{longtable}
\normalsize
(ii)  $Q=1$.
The first table gives  examples where $\mathcal{H}_1\mathcal{H}_2$ and $\mathcal{H}_3\mathcal{H}_4$ are
  principal in $\kk$, when the second shows that $\mathcal{H}_0$,  $\mathcal{H}_1$ and
 $\mathcal{H}_3$ capitulate in $\k$.
\begin{center}
\begin{tabular}{| c | c | c | c | c | }
\hline
$d = 2.p_1.p_2$ & $N(\varepsilon_d )$ & $Q$ & $\mathcal{H}_1\mathcal{H}_2$ in $\kk$ & $\mathcal{H}_3\mathcal{H}_4$ in $\kk$\\
\hline
$890=2.5.89$ & $1$ & $1$ & $[0, 0, 0]~$ & $[0, 0, 0]~$\\
\hline
$1802=2.53.17$ & $1$ & $1$ & $[0, 0, 0]~$ & $[0, 0, 0]~$\\
\hline
$2938=2.13.113$ & $1$ & $1$ & $[0, 0, 0]~$ & $[0, 0, 0]~$\\
\hline
\end{tabular}
\end{center}
\begin{longtable}{| c | c | c | c |  }
\hline
$d = 2.p_1.p_2$ &  $\mathcal{H}_0$ in $\k$ & $\mathcal{H}_1$ in $\k$ & $\mathcal{H}_3$ in $\k$\\
\hline
\endfirsthead
\hline
$d = 2.p_1.p_2$ & $\mathcal{H}_0$ in $\k$ & $\mathcal{H}_1$ in $\k$ & $\mathcal{H}_3$ in $\k$\\
\hline
\endhead
$890=2.5.89$ & $[0, 0, 0, 0, 0]~$ & $[0, 0, 0, 0, 0]~$ & $[0, 0, 0, 0, 0]~$\\
\hline
$1802=2.53.17$  & $[0, 0, 0, 0, 0]~$ & $[0, 0, 0, 0, 0]~$ & $[0, 0, 0, 0, 0]~$\\
\hline
$2938=2.13.113$ &  $[0, 0, 0, 0, 0]~$ & $[0, 0, 0, 0, 0]~$ & $[0, 0, 0, 0, 0]~$ \\
\hline
\end{longtable}
\end{exams}
\section{\textbf{Application}}
Assume $\mathbf{C}_{\mathds{\kk},2}$ is of type $(2 ,2, 2)$;  which occurs if and only if  $p_1\equiv p_2\equiv1\pmod4$ are primes
and at least two  elements of the set \{$(\frac{p_1}{p_2})$, $(\frac{2}{p_1})$, $(\frac{2}{p_2})$\} are equal
to $-1$ (see \cite{AzTa-08}). Under these assumptions the norm of the fundamental unit of
$\QQ(\sqrt{2p_1p_2})$ is equal to $-1$ and  $\mathbf{C}_{\mathds{k},2}=Am_s(\kk/\QQ(i))=\langle [\mathcal{H}_0], [\mathcal{H}_1], [\mathcal{H}_2]\rangle$ (see \cite{AZT12-2}).
\begin{them}\label{43}
Let  $\kk=\QQ(\sqrt{2p_1p_2}, i)$. Put
 $\KK_1=\kk(\sqrt {p_1})$, $\KK_2=\kk(\sqrt {p_2})$,
 then exactly four classes
 capitulate in $\KK_1$,  $\KK_2$ and   $kerJ_{\KK_1}=\langle[\mathcal{H}_1], [\mathcal{H}_2]\rangle$,
 $kerJ_{\KK_2}=\langle[\mathcal{H}_0\mathcal{H}_1], [\mathcal{H}_0\mathcal{H}_2]\rangle$.
\end{them}
\begin{proof}
Since  $N(\varepsilon_3)=-1$, so from Theorems  \ref{31},  \ref{35}  we get
$|kerJ_{\KK_i}|=4$, where  $i\in\{1,2\}$.  Corollaries  \ref{32},  \ref{36}  imply that
 $kerJ_{\KK_1}=\langle[\mathcal{H}_1], [\mathcal{H}_2]\rangle$ and $kerJ_{\KK_2}=\langle[\mathcal{H}_3],
 [\mathcal{H}_4]\rangle$, but as the norm of  $\varepsilon_3$ is equal to $-1$, then
 $\mathcal{H}_0\mathcal{H}_1\mathcal{H}_3$ and $\mathcal{H}_0\mathcal{H}_2\mathcal{H}_4$ or
 $\mathcal{H}_0\mathcal{H}_2\mathcal{H}_3$ and $\mathcal{H}_0\mathcal{H}_1\mathcal{H}_4$ are principal
 in  $\kk$, therefore  $kerJ_{\KK_2}=\langle[\mathcal{H}_0\mathcal{H}_1], [\mathcal{H}_0\mathcal{H}_2]\rangle$.
\end{proof}
\begin{them}\label{44}
Let $\KK_3=\kk(\sqrt {2})=\QQ(\sqrt{2},
\sqrt{p_1p_2}, i)$, then
 \begin{enumerate}[\upshape\indent(1)]
   \item If  $q(\KK_3^+/\QQ)=1$, then four classes
   capitulate in $\KK_3$ and \\$kerJ_{\KK_3}=\langle[\mathcal{H}_0],[\mathcal{H}_1\mathcal{H}_2]\rangle$.
   \item If  $q(\KK_3^+/\QQ)=2$, then  $kerJ_{\KK_3}=\langle[\mathcal{H}_0]\rangle$.
 \end{enumerate}
\end{them}
\begin{proof}
The results are deductions from  Corollary \ref{38}.
\end{proof}
From Theorems \ref{43},  \ref{44} we deduce the following result.
\begin{coro}
 $kerJ_{\k}=\mathbf{C}_{\mathds{k},2}=Am_s(\kk/\QQ(i))$.
\end{coro}

\end{document}